\documentclass[11pt]{article}

\usepackage{url}
\usepackage[usenames,dvipsnames]{color}
\usepackage{fullpage,graphicx,amssymb,amsmath}
\usepackage{cleveref}
\graphicspath{{images/}}
\usepackage{enumerate}
\usepackage{multicol}
\usepackage{multirow}
\usepackage{mathtools}
\usepackage{caption}
\usepackage{subcaption}
\usepackage{ltablex}
\usepackage{relsize}
\usepackage{booktabs}
\usepackage{floatrow}
\usepackage{tabularx}
\usepackage[authoryear, semicolon, sort&compress]{natbib}
\newfloatcommand{capbtabbox}{table}[][\FBwidth]
\renewcommand{\v}[1]{\ensuremath{\mathbf{#1}}}
\def\st{{\rm s.t.}}
\def\-{\raisebox{.75pt}{-}}
\newcommand{\E}{\mathbb{E}}  
\newcommand{\R}{\mathbb{R}}  
\newcommand{\vgamma}{\boldsymbol{\gamma}}
\newcommand{\vdelta}{\boldsymbol{\delta}}
\newcommand{\vlambda}{\boldsymbol{\lambda}}
\newcommand{\vmu}{\boldsymbol{\mu}}
\newcommand{\vxi}{\boldsymbol{\xi}}
\newcommand{\vI}{\boldsymbol{I}}
\newcommand{\vd}{\boldsymbol{d}}
\newcommand{\vV}{\boldsymbol{V}}
\newcommand{\vS}{\boldsymbol{S}}
\newcommand{\vX}{\boldsymbol{X}}
\newcommand{\vY}{\boldsymbol{Y}}
\newcommand{\vh}{\boldsymbol{h}}
\newcommand{\vu}{\boldsymbol{u}}
\newcommand{\vv}{\boldsymbol{v}}
\newcommand{\vA}{\boldsymbol{A}}
\newcommand{\vb}{\boldsymbol{b}}
\newcommand{\vW}{\boldsymbol{W}}
\newcommand{\vw}{\boldsymbol{w}}
\newcommand{\ve}{\boldsymbol{e}}
\newcommand{\vy}{\boldsymbol{y}}
\newcommand{\vz}{\boldsymbol{z}}
\newcommand{\T}{\mathcal{T}} 
\newcommand{\I}{\mathcal{I}} 
\newcommand{\J}{\mathcal{J}}

\begin{document}
\title{Hybrid Strategies using Linear and Piecewise-Linear Decision Rules for Multistage Adaptive Linear Optimization}
\author{Said Rahal$^1$, Dimitri J. Papageorgiou$^{2,\bullet}$, and Zukui Li$^{1,*}$\\
{\small $^1$ Department of Chemical and Materials Engineering}\\
{\small University of Alberta}\\
{\small 9211- 116 St, Edmonton, AB T6G 1H9 CA}\\
{\small $^2$ Corporate Strategic Research}\\
{\small ExxonMobil Research and Engineering Company}\\
{\small 1545 Route 22 East, Annandale, NJ 08801}\\
{\small $^*$\texttt{zukui@ualberta.ca}, $^\bullet$\texttt{dimitri.j.papageorgiou@exxonmobil.com}} \\
~\\
}
\date{}
\maketitle
\begin{abstract}

Decision rules offer a rich and tractable framework for solving certain classes of multistage adaptive optimization problems. Recent literature has shown the promise of using linear and nonlinear decision rules in which wait-and-see decisions are represented as functions, whose parameters are decision variables to be optimized, of the underlying uncertain parameters.  Despite this growing success, solving real-world stochastic optimization problems can become computationally prohibitive when using nonlinear decision rules, and in some cases, linear ones. Consequently, decision rules that offer a competitive trade-off between solution quality and computational time become more attractive. Whereas the extant research has always used homogeneous decision rules, the major contribution of this paper is a computational exploration of hybrid decision rules. We first verify empirically that having higher uncertainty resolution or more linear pieces in early stages is more significant than having it in late stages in terms of solution quality. Then we conduct a comprehensive computational study for non-increasing (i.e., higher uncertainty resolution in early stages) and non-decreasing (i.e., higher uncertainty resolution in late stages) hybrid decision rules to illustrate the trade-off between solution quality and computational cost. We also demonstrate a case where a linear decision rule is superior to a piecewise-linear decision rule within a simulator environment, which supports the need to assess the quality of decision rules obtained from a look-ahead model within a simulator rather than just using the look-ahead model's objective function value.


\end{abstract}

\newpage

\section{Introduction}
        In recent years, multistage adaptive optimization has received growing interest as a tool to address parameter uncertainty in decision making problems. In particular, multistage adaptive optimization can be cast as a sequential decision making problem under uncertainty; see \cite{georghiou2015generalized}. Instead of taking all decisions at once, without any previous knowledge of the progressively revealed uncertain parameters, the decision maker first implements a static decision $\v{x}_1$ which is known to be independent of the future uncertainty values. Afterwards, the decision maker waits for the gradual unfolding of uncertainty to implement the optimal recourse decision. In this sequence, the first uncertain parameter $\vxi_2$ is revealed\footnote{The time subscript refers to the stage by which the information of a given variable is available. The first observed uncertainty $\vxi_2$ is first available in stage 2.}, followed by a recourse decision $\v{x}_2$. Practically, the recourse decision is a function of the realized uncertainty $\v{x}_2 \equiv \v{x}_2(\vxi_2)$.
After which, the sequence of alternating observations and recourse decisions unfolds over $T$ stages, where in each stage $t \in \T_{-1}= \{2,\ldots,T\}$, the decision maker observes an uncertain parameter $\vxi_t$ and selects a recourse decision $\v{x}_t(\vxi_{[t]})$,  a decision which depends on the whole history of past observations $\vxi_{[t]} = (\vxi_2,\dots,\vxi_{t})$,  but not on any future observations $\vxi_{t+1},\dots,\vxi_T$. According to \cite{shapiro2009lectures}, a general formulation of a multistage adaptive optimization problem is given as
\begin{subequations} \label{model:generic}
	\begin{alignat}{4}
	\min_{\v{x}_t(\cdot)}~~& \v{c}_1^{\top}\v{x}_1 + \rho\left[ \sum_{t=2}^T\v{c}_t(\vxi_{[t]})^{\top} \v{x}_t(\vxi_{[t]}) \right] &&  \\
	\st~~ & \v{A}_1\v{x}_1 \geq \v{b}_1 &&  \\
	      & \sum_{s=2}^t \v{A}_{s}(\vxi_{[s]})  \v{x}_s(\vxi_{[s]})  \geq \v{b}_t(\vxi_{[t]}) & & \qquad \forall \vxi \in \Xi,\ t \in \T_{-1}  
	\end{alignat}
\end{subequations}
where $\vA_1$, $\vb_1$, $\v{c}_1$ are the first-stage static parameters and $\v{x}_1$ is the first-stage static (or \textit{here and now}) decisions. Similarly, $A_t(\vxi_{[t]})$, $b_t(\vxi_{[t]})$, $c_t(\vxi_{[t]})$ are uncertain recourse matrices, right-hand side vectors, cost coefficients and the functionals $\v{x}_t(\vxi_{[t]})$ are the adaptive (or recourse or \textit{wait and see}) decisions. The possible realizations of $\vxi$ are defined by the underlying uncertainty set $\Xi$. The functional $\rho$ is a coherent risk measure \citep{Rockafellar2007}. 

There are numerous paradigms to tackle sequential decision making problems. Arguably the three most popular within the operations research community are stochastic programming \citep{shapiro2009lectures}, adaptive robust optimization \citep{ben2009robust}, and stochastic dynamic programming \citep{powell2011_adp_book}. The first two grew out of a mathematical programming tradition, while the third has roots in control theory and reinforcement learning. While all paradigms aim to compute optimal decision policies, they do so based on assumptions about the underlying uncertainty required to define $\rho$. Stochastic programming typically requires some underlying distributional assumptions about the uncertain parameters to determine the optimal expected objective function value, while robust optimization only assumes that the uncertain parameters belong to a known uncertainty set to compute the worst-case scenario. Meanwhile, stochastic dynamic programming commonly relies on the solution of Bellman's equation, or approximations thereof, to generate policies.  Today, one could argue that the lines between the three domains are becoming even more blurred as a cross-pollination of ideas continues to flourish.

Of the numerous algorithmic advances devised to construct optimal policies, two of the most prominent are scenario- and decision rule (DR)-based methods. The former is typically used in the context of stochastic programming where a set of discrete scenarios representing the uncertainty set is used to compute the expected objective function value. Scenarios are often represented via a scenario tree whose size increases exponentially with the size of the decision making sequence. Such increase induces prohibitive computational overhead which is the main limitation of scenario-based stochastic programming methods, and is known as curse of dimensionality.

On the other hand, decision rule-based methods do not suffer from the same curse of dimensionality. First introduced by \cite{garstka1974decision}, if not earlier, the significance of the approach was not fully realized until 2004 when \cite{ben2004adjustable} demonstrated that the best linear decision rule (LDR) for robust and stochastic optimization can be solved in polynomial time. In their framework, \cite{ben2004adjustable} defined uncertainty-independent decision variables as static, while uncertainty-dependent adaptive functional decisions are defined as linear functions, or rules, of the uncertain parameters. The main advantage of LDRs is that, under certain convexity assumptions of the underlying uncertainty set, they give rise to tractable robust counterparts that are efficiently solvable by today's ever-improving optimization engines. However, this modelling feature comes at the expense of low solution quality. \cite{lappas2016multi} provide an instructive overview of implementing LDRs in process scheduling problems in the context of robust optimization.

Variants of nonlinear decision rules which improve the solution quality of the adaptive functional decisions have been proposed in the literature. \cite{chen2008linear} introduced deflected and segregated decision rules used for stochastic programming problems with semi-complete and general recourse. \cite{see2010robust} implemented a truncated LDR for their inventory problem which was proven to do better than the LDR policy. Both decision rules are generalized in \cite{goh2010distributionally} as bi-deflected LDRs. \cite{ben2011immunizing}  restricted the functional decisions into quadratic decision rules in the context of robust optimization. The obtained tractable counterpart, under an ellipsoidal uncertainty set, is a second order cone programming problem. Further, \cite{bertsimas2011hierarchy} proposed polynomial decision rules in multistage robust dynamic problems. The increase in solution quality comes at the cost of increased computational complexity; the tractable counterpart, under an intersection of convex uncertainty sets, is a semidefinite programming problem. Polynomial decision rules were later refined by \cite{bampou2011scenario} in the context of stochastic programming problems.

A specific class of nonlinear DRs is the piecewise linear decision rule (PLDR). It improves the flexibility of a DR by having multiple slopes (i.e., decision variables) while inheriting the modelling features of an LDR due to its linear nature.  \cite{chen2009uncertain} proposed an extended linear decision rule using an extended uncertainty set defined via the positive and negative perturbations of the original uncertainty set. The decision rule is equivalent to a PLDR with two linear pieces. Later, \cite{georghiou2015generalized} introduced the concept of  generalized decision rules via liftings. The key property is the one to one correspondence between the linear decision rules in the lifted problem and a family of nonlinear decision rules in the original problem. Hence, the modelling features of LDRs in the lifted space are exploited, while still exhibiting the flexibility of nonlinear decision rules in the original space. A PLDR is an example of lifted decision rules; they are derived for any number of breakpoints or linear pieces in \cite{georghiou2015generalized}. Recently, \cite{ben2016tractable} proposed a novel piecewise linear decision rule for linear dynamic robust problems. The framework is based on approximating the uncertainty set with a simplex where constructing the PLDR, with exponential number of pieces, can be performed efficiently. For a recent comprehensive survey of decision rules, see \cite{YANIKOGLU2018}. 

Since our focus is on hybridizing the lifting in PLDRs, it is worth calling attention to several noteworthy applications where PLDRs have been employed to demonstrate their growing popularity.  In \cite{munoz2014piecewise}, PLDRs were used to approximate recourse decisions when dispatching electric power given random power supply and consumption. \cite{GAUVIN2017} evaluated various LDRs and PLDRs for managing reservoirs in Canada for electric power generation. \cite{braathen2013hydropower} compared scenario-, LDR- and PLDR-based approaches in optimizing the hydropower bidding process for Nordic producers. \cite{Limeng2015IterativeLDR} extended the application of LDRs and segregated DRs to nonlinear concave objective functions (in a maximization problem) for optimal reservoir operation.  \cite{beuchat2016performance} integrated power dispatch  and reserve models and illustrated that the least flexible PLDR (specifically, two pieces for each uncertainty dimension) provides substantial performance improvement with respect to an LDR. Further, in the control community,  \cite{jin2018segregated} implemented segregated LDRs to solve multistage stochastic control problem with linear dynamics and quadratic cost.  \cite{zhang2015convex} considered linear chance constrained model predictive control problems subject to additive disturbance. For their problem, both the randomization approach (i.e., scenario-based approach) and the PLDR-based method were found to be computationally expensive. As an alternative, they implemented a combination of the two approaches which allowed them to exploit the flexibility of a PLDR. 

The choice of which decision rule to implement is governed by two competing objectives: high solution quality and low computational cost. The latter objective favours LDRs while the former is portrayed more by nonlinear DRs. To the best of our knowledge, decision rule-based methods found in the literature implement homogeneous decision rules That is, the same form of linear or nonlinear decision rule is applied to every adaptive decision variable (in every stage) that is chosen to be represented with a decision rule. This observation motivates our work to investigate hybrid decision rules combining the salient features of both types. We will limit our study of nonlinear DRs to PLDRs only. For example, lifting the uncertainty set for near-term stages (e.g., stages $2-4$) and keeping the original uncertainty for the remaining stages (e.g., stages $5 -T$) gives increased flexibility in the types of permissible decisions/actions in the immediate future while giving limited recourse actions in subsequent long-term decisions. It is somewhat analogous to an approach used in scenario-based stochastic programming in which a scenario tree containing many branches per node (i.e., higher uncertainty resolution) in early stages and relatively few branches per node (i.e., lower uncertainty resolution) in later stages can be attractive (\cite{bakkehaug2014stochastic}, \cite{arslan2017bulk}).

The contributions of this paper are:
\begin{enumerate}
\item Similar to what has been illustrated in scenario-based stochastic programming methods, we empirically show that ``it is more important to model the uncertainty of the near future with more details than it is for the later stages'' (\cite{bakkehaug2014stochastic}, p.73). We demonstrate this result for the first time using decision rules where having higher uncertainty resolution or more linear pieces in early stages improves the flexibility of a policy more than having it in late stages. This observation is shown via an empirical sensitivity analysis for two computational settings with two different planning horizons, which motivates the design of PLDRs with axial segmentation using hybrid combination of liftings (i.e., HDRs). It is important to credit \cite{georghiou2015generalized} with conceiving the idea of exploiting the modularity of decision rules in stochastic programming.  At the same time, these authors did not pursue this modular design with detailed computational experiments. Our computational study attempts to fill this gap.
\item We perform a comprehensive computational study on the design of hybrid decision rules with a non-increasing (i.e., higher uncertainty resolution or more linear pieces in early stages) and non-decreasing (i.e., higher uncertainty resolution or more linear pieces in late stages) lifting strategies. We demonstrate empirically, using marginal distribution plots, that (1) non-increasing HDRs are always more flexible than non-decreasing HDRs in terms of solution quality, and (2) the computational benefits of non-increasing HDRs as competitive candidates for the trade-off between solution quality and computational time is best manifested using a well designed lifting strategy. As in any design problem, we show that a poorly designed non-increasing HDR may lose the aforementioned computational benefits. We systematically illustrate the impact of the lifting strategy of a non-increasing HDR on the trade-off between solution quality and computational time.
\item We show that a linear decision rule can be \textit{superior} to a  piecewise linear decision rule with axial segmentation and a single breakpoint within a simulator environment. This observation is counter-intuitive as PLDRs are more flexible than LDRs. The reason for this peculiar behaviour is due to (1) the presence of \textit{mutually exclusive} state variables and (2) the robust nature of the stochastic counterparts. Our case study reveals that it is crucial to evaluate policies within a simulation environment to obtain an impartial assessment of how various policies perform in practice. \cite{powell2014clearing} also advocates this point. 
\end{enumerate}

The rest of the paper is organized as follows. In section 2, we illustrate the derivation of the tractable stochastic counterparts using LDRs and PLDRs for a multistage stochastic newsvendor problem. In section 3, we showcase a specific newsvendor problem setting where an LDR is superior to a PLDR with a single breakpoint. In section 4, we illustrate that having higher uncertainty resolution in early stages of a hybrid decision rule is more important than in late stages. We demonstrate the computational benefits acquired by non-increasing hybrid decision rules through a set of comprehensive computational experiments. In section 5, we conclude the paper and offer future research directions. 

\section{Linear vs. piecewise-linear decision rules}
Throughout our study, we assume that the cost coefficients and recourse matrices are fixed and the risk measure $\rho$ in Model \eqref{model:generic} is the expectation functional. The resulting multistage stochastic adaptive problem has the following simplified form
\begin{subequations} \label{model:generic_simplified}
	\begin{alignat}{4}
	\min_{\v{x}_t(\cdot)}~~& \v{c}_1^{\top}\v{x}_1 + \E\left[ \sum_{t=2}^T\v{c}_t^{\top} \v{x}_t(\vxi_{[t]}) \right] &&  \\
	\st~~ & \v{A}_1\v{x}_1 \geq \v{b}_1 &&  \\\label{model:generic_simplified_semiinf}
	      & \sum_{s=2}^t \v{A}_{s}\v{x}_s(\vxi_{[s]})  \geq \v{b}_t(\vxi_{[t]}) & & \qquad \forall \vxi \in \Xi,\ t \in \T_{-1}   
	\end{alignat}
\end{subequations}
Formulation \eqref{model:generic_simplified} is computationally intractable due to the presence of semi-infinite constraints. Decision rule based-methods circumvent this intractability by defining the adaptive decisions $\v{x}_t(\vxi_{[t]})$ as a specific function or rule of the uncertain parameters. The simplest rule is an LDR where  $\v{x}_t(\vxi_{[t]})$ is defined as 
\begin{equation*}
  \v{x}_t(\vxi_{[t]}) = \v{x}_{t}^0 + \sum_{s=2}^t \v{X}_{s}^{1 } \vxi_s \quad \ \forall \vxi \in \Xi,\ t \in \T_{-1}
\end{equation*}
where $\v{x}_t^0$ and $\v{X}_{t}^1$ are the intercepts and slopes, respectively. As is typically done, we assume that $\v{b}_t(\vxi_{[t]})$ is a linear/affine function of the uncertain parameters $ \v{b}_t(\vxi_{[t]}) = \v{b}^0_t + \sum_{s=2}^t \v{B}_{s}^{1}\vxi_s$, where $B_s^1$ define this linear dependence. Substituting for $\v{x}_t(\vxi_{[t]})$ and $\v{b}_t(\vxi_{[t]})$ in Model \eqref{model:generic_simplified}, we obtain
\begin{subequations}
	\begin{alignat}{4}
	\min_{\v{x}_t(\cdot)}~~& \v{c}_1^{\top}\v{x}_1 +  \sum_{t=2}^T\v{c}_t^{\top} \left(\v{x}_{t}^0 + \sum_{s=2}^t \v{X}_{s}^{1} \E[\vxi_s]\right) &&  \\
	\st~~ & \v{A}_1\v{x}_1 \geq \v{b}_1 &&  \\
	      & \sum_{s=2}^t \v{A}_{s}\left(\v{x}_{s}^0 + \sum_{p=2}^s \v{X}_{p}^{1} \vxi_p \right) \geq \v{b}^0_t + \sum_{s=2}^t \v{B}_{s}^{1}\vxi_s & & \qquad \forall \vxi \in \Xi,\ t \in \T_{-1}  
	\end{alignat}
\end{subequations}
Under certain convexity assumption of the uncertainty set $\Xi$, the tractable counterpart is derived by exploiting the strong duality property of convex optimization problems.

Extending the idea above, \cite{georghiou2015generalized} introduced PLDRs with axial segmentation. The decision rule is still linear in the lifted uncertainty space, however it corresponds to a PLDR in the original uncertainty space. In this regard, $\v{x}_t(\vxi^\prime_{[t]})$ is defined as 
  $$
  \v{x}_t(\vxi^\prime_{[t]}) = \v{x}_{t}^{\prime 0} + \sum_{s=2}^t \v{X}_{s}^{\prime 1} \vxi^\prime_s \quad \ \forall \vxi^\prime \in \Xi^\prime,\ t \in \T_{-1}
  $$
where $\v{x}_t^{\prime 0}$ and $\v{X}_{t}^{\prime 1}$ are the intercepts and slopes in the lifted space, respectively. The lifted uncertainty set $\Xi^\prime$ defines the possible realizations of $\vxi^\prime$.

\begin{figure}[H]
\centering
\begin{subfigure}{0.465\textwidth}
\centering
\includegraphics[scale=0.39]{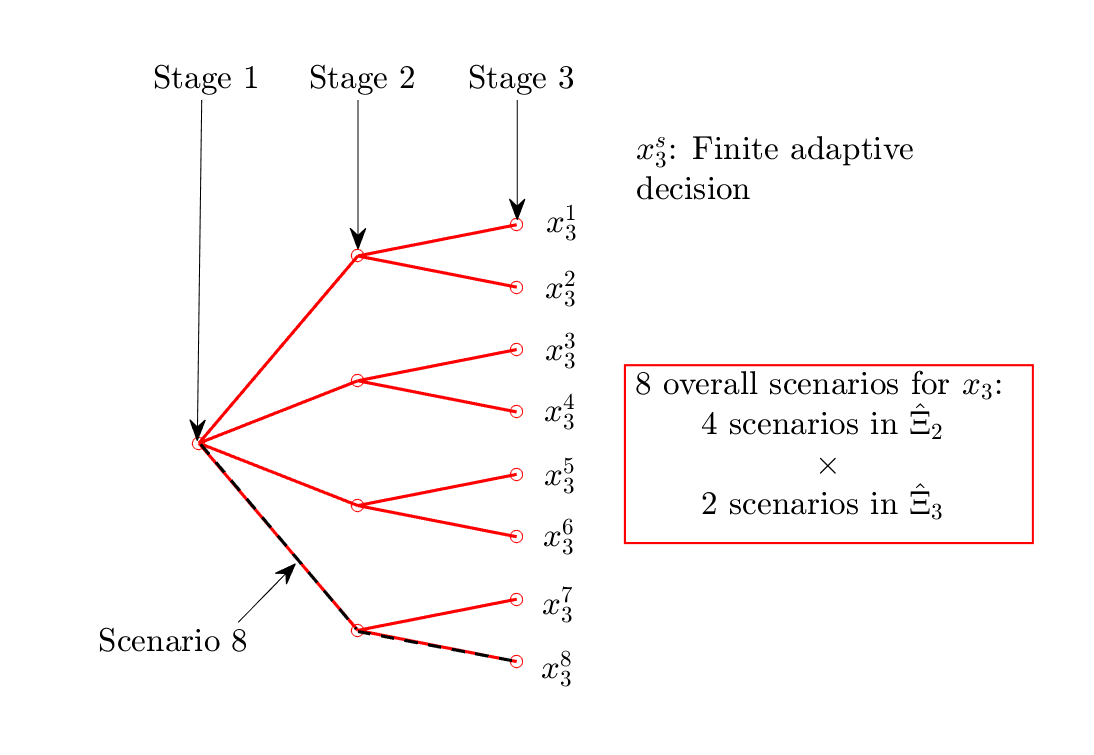}
\caption{Three-stage scenario tree}
\end{subfigure}
\begin{subfigure}{0.525\textwidth}
\centering
\includegraphics[scale=0.39]{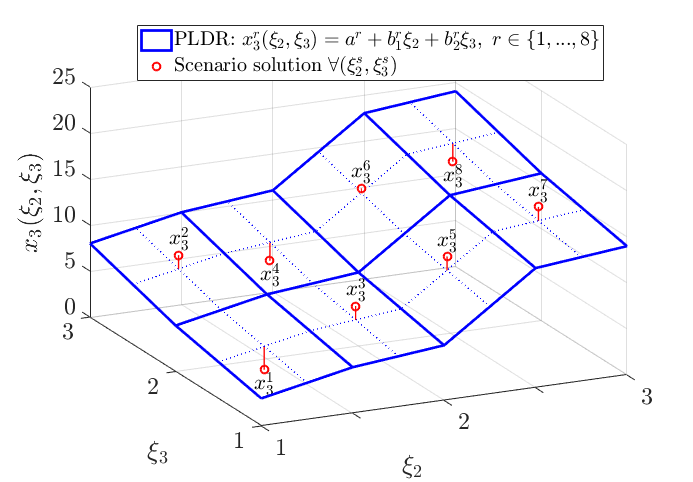}
\caption{Scenario- and DR-based recourse decisions in stage $3$}
\end{subfigure}
\caption{Scenario- and decision rule-based recourse decisions in stage $3$ of a multistage adaptive optimization problem where $\xi_2,\ \xi_3 \in \R$. In the scenario-based approach [Fig. (a)], eight recourse decisions $x_3^s$ are generated due to four and two scenarios in stages 2 and 3, respectively. Similarly, in the DR-based approach [Fig. (b)], eight recourse functions $x_3^r(\xi_2,\xi_3)$ are obtained by lifting $\xi_2$ and $\xi_3$ to four and two elements, respectively. Figure (b) illustrates the possible deviation between the recourse decisions $x_3^s$ and the recourse functions $x_3^r(\xi_2,\xi_3)$, which reflects the different approximations resulting from the two approaches. However, the comparable impact of scenarios on $x_3^s$ and of lifting on $x_3^r(\xi_2,\xi_3)$ motivate our work to verify empirical evidence found in scenario-based stochastic programming methods in the context of decision rule-based methods, where modelling uncertainty with higher resolution in early stages has been shown to be more attractive. }
\label{fig:analogy_scenario_decision_rule}
\end{figure}

Interestingly, some similarities can be drawn between PLDR-based methods and scenario-based stochastic programming methods. In the latter approach, the uncertainty set is approximated by a discrete set of scenarios for which a discrete set of optimal recourse decisions is computed. However, a common feature of both approaches is that the solution quality of the recourse decisions increases as the granularity/resolution of the uncertainty increases (be it more discrete scenarios or more lifted elements). To this end, Figure \ref{fig:analogy_scenario_decision_rule} attempts to graphically contrast the two approaches in a 3-stage example by illustrating the solution of the recourse decision in stage 3 where the uncertain parameters $\xi_2$ and $\xi_3$ are one-dimensional. While scenario-based stochastic programming methods consider discrete uncertainty sets $\hat{\Xi}_2$ and $\hat{\Xi}_3$  of scenarios and determine a recourse action for each scenario, decision rule-based methods provide an infinite number of recourse actions - one for every realization of the uncertain parameter in $\Xi_2$ and $\Xi_3$. Still, the impact of the number of scenarios and lifted elements on the recourse decision $x_3$ is comparable. For example, implementing four scenarios in stage 2 and two scenarios in stage 3 generate eight possible finite decisions $x_3^s$, where $x_t^s$ denotes a recourse decision from a scenario-based method in stage $t$ for scenario $s$.  Likewise, introducing four lifted elements in stage 2 and two lifted elements in stage 3 generates eight linear decision functions $x^r_3(\xi_2,\xi_3)$, where $x_t^r$ denotes a recourse function  from a PLDR-based method in stage $t$ and within subspace $r$ of the lifted uncertainty space. 
Throughout this paper, we make the following assumptions:
\begin{itemize}
\item \textbf{Assumption 1:} Piecewise linear decision rules and hybrid decision rules are constructed via lifting with \textit{axial segmentation} as described in \cite{georghiou2015generalized}. We do not address lifting with generalized segmentation.
\item \textbf{Assumption 2:} The set of potential breakpoints to construct PLDRs is given. The search for an optimal set of breakpoints in each stage is still an open question, but is out of scope for this computational study.
\item \textbf{Assumption 3:} The set of breakpoints implemented in a PLDR and for a specific resolution in an HDR is \textit{the same in all stages}. Clearly, this may not be optimal, but it allows us to perform comprehensive computational experiments.
\end{itemize}

In this section, we introduce key concepts about LDR and PLDR methods by way of example. 
We first introduce a multistage stochastic newsvendor problem to illustrate how decision rules can be applied.  We then derive the stochastic counterparts using an LDR and a PLDR. Finally, we compare LDRs and PLDRs  and illustrate the improvement induced by the additional flexibility of PLDRs.

\subsection{Illustrative example: multistage stochastic newsvendor problem}
\begin{figure}[H]
\begin{center}
\includegraphics[scale=0.6]{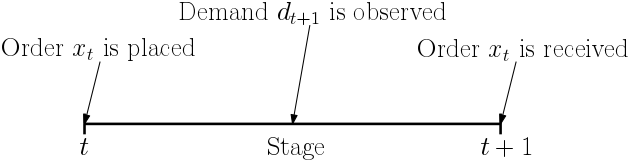}
\caption{Chronological sequence of placing an order $x_t$, observing the demand $d_{t+1}$, and receiving the order in a multistage stochastic newsvendor problem. The order is placed before observing the demand, however it is being recieved afterwards.} \label{fig:time_sequence_vars_param_multistage_inv}
\end{center}
\end{figure}

In a multistage newsvendor problem, a seller has to satisfy the demand of a perishable good at the minimal total cost. At each stage $t\in \T_{-T}$, an order $x_t$ is placed to satisfy the demand $d_{t+1}$ and is first available for selling at the next stage $t+1$. The cumulative difference between $x_t$ and $d_{t+1}$ defines the inventory and backlog: positive value indicates inventory amounts, while negative value indicates backlog amounts.
Figure \ref{fig:time_sequence_vars_param_multistage_inv} illustrates the chronological order of $x_t$ with respect to $d_{t+1}$; the latter is observed after the order $x_t$ is placed, but before it is received.

The deterministic formulation of a multistage newsvendor problem is
\begin{subequations} \label{model:invprob_det_multist}
	\begin{alignat}{4} \label{model:inv_det_obj}
	 \min_{x_t,I_t} ~~  & \sum_{t=1}^{T-1}C_tx_t + \sum_{t=2}^{T}(H_t (I_{t})^+  + B_t (-I_{t})^+) && \\\label{model:inv_det_inv}
	 \st ~~ & I_{t} = I_{t-1}+ x_{t-1} - d_{t}  && \qquad \forall t \in \T_{-1}     \\\label{model:inv_det_order}
            &  0 \leq x_t \leq U^x       && \qquad \forall t \in \T_{-T}
	\end{alignat}
\end{subequations}
where $I_1$ is the initial inventory and $[\cdot]^+ = \max(0,\cdot)$. The objective function includes ordering, holding and backlogging costs. Ordering too much results in high ordering and additional holding costs, whereas ordering too little incurs expensive backlogging costs. The inventory balance is governed by  eq. \eqref{model:inv_det_inv}. At last, eq. \eqref{model:inv_det_order} defines the lower and upper bounds of the ordering amount.

After introducing auxiliary variables, Model \eqref{model:invprob_det_multist} can be reformulated as a linear programming problem
\begin{subequations} \label{model:news_det_multist_2}
	\begin{alignat}{4}
	 \min_{\substack{x_t,I_t \\s_t^+,s_t^-}}~~ & \sum_{t=1}^{T-1}C_tx_t + \sum_{t =2}^{T}(H_t s_t^+ + B_t s_t^-) && \\\label{model:inv_det_inv_2}
	 \st ~~ & I_{t} = I_{t-1}+ x_{t-1} - d_{t}  &&  \qquad \forall t \in \T_{-1}   \\ \label{model:inv_det_inv_3}
	        & s_t^+ \geq I_{t}    && \qquad \forall t \in \T_{-1}   \\ \label{model:inv_det_inv_4}
	        & s_t^- \geq -I_{t} && \qquad \forall t \in \T_{-1}   \\
	        & 0 \leq x_t \leq U^x && \qquad \forall t \in \T_{-T}\\ 
            & s_t^+,s_t^- \ge 0 && \qquad \forall t \in \T_{-1}  
	\end{alignat}
\end{subequations}
where $s_t^+$ and $s_t^-$ are the inventory and backlog amounts, respectively.

In real world applications, the demand may be uncertain. Solving a deterministic model using expected demand will likely yield a suboptimal or potentially infeasible policy. Modelling the uncertainty as the demand parameter itself, a multistage adaptive stochastic newsvendor problem is formulated as
\begin{subequations} \label{model:inv_SP_multist_gen_eq}
	\begin{alignat}{4}
	\label{model:inv_SP_multist_gen_objfct}
	 \min_{\substack{x_t(\cdot),I_t(\cdot) \\s_t^+(\cdot),s_t^-(\cdot)}}~~ & \E \left[\sum_{t=1}^{T-1}C_tx_t(\vd_{[t]}) + \sum_{t=2}^{T}(H_t s_t^+(\vd_{[t]}) + B_t s_t^-(\vd_{[t]}))\right] &&  \\\label{model:inv_SP_multist_gen_start}
     \st~~  &  I_t(\vd_{[t]}) = I_{t-1}(\vd_{[t-1]}) +x_{t-1}(\vd_{[t-1]})-d_{t} && \quad \forall \vd \in \Xi,\ t \in \T_{-1} \\	
	& s_t^+(\vd_{[t]}) \geq I_t(\vd_{[t]}) && \quad \forall \vd \in \Xi,\ t \in \T_{-1} \\
	& s_t^-(\vd_{[t]}) \geq - I_t(\vd_{[t]}) && \quad \forall \vd \in \Xi,\ t \in \T_{-1} \\
	& 0 \leq x_t(\vd_{[t]}) \leq U^x  && \quad \forall \vd \in \Xi,\ t \in \T_{-T}          \\ \label{model:inv_SP_multist_gen_end}
	& s_t^+(\vd_{[t]}),\ s_t^-(\vd_{[t]}) \ge 0 && \quad \forall \vd \in \Xi,\ t \in \T_{-1}
	\end{alignat}
\end{subequations}
where $\vd_{[t]}=[d_2,\cdots,d_{t}]$,  $x_1(\vd_{[1]}) \equiv x_1$ is the first-stage ordering decision, $I_1(\vd_{[1]}) \equiv I_1$, and $\Xi$ is the uncertainty set for demand. The expectation is computed with respect to the distribution of $\vd$.

\subsection{Linear adaptive stochastic counterpart of the newsvendor problem}
The adaptive functional decisions $x_t(\vd_{[t]})$ are defined in terms of the past demand $\vd_{[t]} = (d_2,\dots,d_t)$, which are observed and known in stage $t$; they cannot be a function of future unrealized demand parameters $(d_{t+1},\dots,d_T)$  \citep{ben2004adjustable}. To make this temporal dependence concrete, an observation matrix $\vV_t \in \Re^{(T-1)\times (T-1)}$ which relates  $\vd$ to $\vd_{[t]}$ is introduced
\begin{equation}
\vd_{[t]} = \left[\begin{array}{ll}
\vI_{t-1} & \mathbf{0}_{(t-1) \times (T-t)} \\
\mathbf{0}_{(T-t) \times (T-1)} & \mathbf{0}_{(T-t) \times (T-t)} \\
\end{array}\right] \vd = \vV_t \vd = [d_2,d_3,\dots,d_{t},0,\cdots,0]^\top \quad \forall t \in \T_{-1}
\end{equation}
where $\vI_{t-1} \in \Re^{(t-1)\times(t-1)}$ is the identity matrix. 

For $x_t(\vd_{[t]})$, $I_t(\vd_{[t]})$, $s_t^+(\vd_{[t]})$ and  $s_t^-(\vd_{[t]})$, consider the following LDRs
\begin{subequations} \label{model:inv_SP_multist_LDR_eq}
	\begin{alignat}{4}
	& x_t (\vd_{[t]}) =&& ~x_t^0 + \vX_t^1\vV_t\vd && \qquad \forall t \in \T_{-T}\\
	& I_t (\vd_{[t]}) =&& ~I_t^0 + \vI_t^1\vV_t\vd && \qquad \forall t \in \T \\
	& s_t^+(\vd_{[t]})=&& ~s_t^{0+} + \vS_t^{1+}\vV_t\vd && \qquad \forall t \in \T_{-1} \\
	& s_t^-(\vd_{[t]})=&& ~s_t^{0-} + \vS_t^{1-}\vV_t\vd && \qquad \forall t \in \T_{-1}
	\end{alignat}
\end{subequations}
where $x_t^0,\ I_t^0,\ s_t^{0+},\ s_t^{0-}$ are the intercepts and $\vX_t^1,\vS_t^{1+},\ \vS_t^{1-},\ \vI_t^1 \in \Re^{1 \times (T-1)}$ are the slopes. We let $x_{1}^0=x_1$, $I_1^0=I_1$ and $\vX_1^1=\vI_1^1=\textbf{0}$. 

Implementing LDRs in the set of constraints of Model \eqref{model:inv_SP_multist_gen_eq}, we obtain
\begin{subequations} \label{model:newsvendor_affine} 
	\begin{alignat}{4}
\label{eq:newsvendor_affine_invBalance}
	 & I_t^0 +\vI_t^1\vV_t\vd= I_{t-1}^0 +\vI_{t-1}^1\vV_{t-1}\vd +x_{t-1}^0 +\vX_{t-1}^1\vV_{t-1}\vd +d_{t} && \quad \forall \vd \in \Xi,\ t \in \T_{-1} \\	\label{eq:newsvendor_affine_storage}
    & s_t^{0+} + \vS_t^{1+}\vV_t\vd \ge I_t^0 + \vI_t^1\vV_t\vd &&  \quad \forall \vd \in \Xi,\ t \in \T_{-1} \\
	& s_t^{0-} + \vS_t^{1-}\vV_t\vd \ge -I_t^0-\vI_t^1\vV_t\vd && \quad \forall \vd \in \Xi,\ t \in \T_{-1} \\
		& x_t^0 + \vX_t^1\vV_t\vd \geq 0      &&  \quad \forall \vd \in \Xi,\ t \in \T_{-T} \\
	& x_t^0 + \vX_t^1\vV_t\vd \leq U^x    &&  \quad \forall \vd \in \Xi,\ t \in \T_{-T} \\
	& s_t^{0+} + \vS_t^{1+}\vV_t\vd \ge 0 &&  \quad \forall \vd \in \Xi,\ t \in \T_{-1} \\
	& s_t^{0-} + \vS_t^{1-}\vV_t\vd \ge 0 &&  \quad \forall \vd \in \Xi,\ t \in \T_{-1} 
	\end{alignat}
\end{subequations} 
We assume that the uncertainty set $\Xi$ is defined by a generic polytope:\begin{equation} \label{eq:multistg_inv_polyhedral_uncertainty_set}
	\Xi := \left\lbrace\vd \in \Re^{T-1}|\ \vW\vd \ge \vh \right\rbrace
\end{equation}
where $\vW \in \Re^{m\times (T-1)}\ \text{and}\ \vh \in \Re^m$. As a mean of example, we will derive the stochastic counterpart for the semi-infinite constraint in eq. \eqref{eq:newsvendor_affine_storage}. First, it is rearranged as follows
	\begin{equation} \label{eq:multistg_inven_const1_equiv}
 	  s_t^{0+} - I_t^0  + \left\{\begin{array}{lll}
	\min\limits_{\vd}  &  \left(\vS_t^{1+}\vV_t-\vI_t^1\vV_t\right) \vd \\
	\st	 & \vW\vd \ge \vh  \end{array}\right\}\geq 0 \qquad \forall t \in \T_{-1}
	\end{equation}
Introducing the $\min_{\vd}(\cdot)$ operator does not affect the solution. When eq. \eqref{eq:multistg_inven_const1_equiv} is satisfied, it follows that eq.  \eqref{eq:newsvendor_affine_storage} is satisfied for all $\vd \in \Xi$. The dual of the inner minimization problem is derived as
    \begin{equation} \label{eq:multistg_inven_const1_dual}
 	 s_t^{0+} - I_t^0  + \left\{\begin{array}{lll}
	\max\limits_{\vu_t}  & \vh^{\top} \vu_t \\
	\st	 & \vW^{\top}\vu_t = \left(\vS_t^{1+}\vV_t-\vI_t^1\vV_t\right)^{\top}\\
	     & \vu_t \in \Re^m_+  \end{array}\right\}\geq 0  \qquad \forall t \in \T_{-1}
	\end{equation}
	where $\vu_t$ is the dual variable. The $\max_{\vu_t}(\cdot)$ operator can be removed without affecting the optimal solution. The final form of the stochastic counterpart of eq. \eqref{eq:newsvendor_affine_storage} is equivalent to
	\begin{subequations} \label{eq:multistg_inven_const1_tractable_counterpart}
	\begin{alignat}{4}
	& s_t^{0+} - I_t^0  + \vh^{\top}\vu_t  \ge 0 && \quad \forall t \in \T_{-1} \\
  & \vW^{\top}\vu_t = \left(\vS_t^{1+}\vV_t-\vI_t^1\vV_t\right)^{\top}  && \quad \forall t \in \T_{-1} \\
  & \vu_t \in \Re^m_+  && \quad \forall t \in \T_{-1}
  \end{alignat}
\end{subequations}	
The stochastic counterparts of the remaining semi-infinite inequality constraints are derived following the same procedure. Equation \eqref{eq:newsvendor_affine_invBalance} is the only semi-infinite equality constraint, and its tractable counterpart is derived by forcing the intercept and the slope of the constraint to be equal to zero.  To illustrate, eq. \eqref{eq:newsvendor_affine_invBalance}  is equivalently rewritten as
\begin{equation} 
I_t^0 -I_{t-1}^0 -x_{t-1}^0 + (\vI_t^1\vV_t - \vI_{t-1}^1\vV_{t-1} - \vX_{t-1}^1\vV_{t-1}+\ve^{\top}_{t-1})\vd = 0 \qquad \forall \vd \in \Xi,\ t \in \T_{-1}
\end{equation}
The vector $\ve_{t} \in \Re^{1\times (T-1)}$ has a value of $1$ at the $t$ index and a value of 0 elsewhere. It is satisfied for all $\vd \in \Xi$ if and only if
\begin{subequations} \label{eq:newsvendor_equality_count}
	\begin{alignat}{4}	
  & I_t^0 -I_{t-1}^0 -x_{t-1}^0  = 0&& \qquad \forall t \in \T_{-1}\\
  & \vI_t^1\vV_t - \vI_{t-1}^1\vV_{t-1} - \vX_{t-1}^1\vV_{t-1}+\ve^{\top}_{t-1} = \v{0} && \qquad \forall t \in \T_{-1}
	\end{alignat}
\end{subequations}
Equation \eqref{eq:newsvendor_equality_count} is considered the stochastic counterpart of eq. \eqref{eq:newsvendor_affine_invBalance}. 

The linear adaptive stochastic counterpart (LASC) of Model \eqref{model:newsvendor_affine} is formulated as
\begin{subequations} \label{model:newsvendor_affine_counterpart} 
	\begin{alignat}{4}
 \label{model:inv_SP_multist_LDR_eq_count_start}
	   & I_t^0 -I_{t-1}^0 -x_{t-1}^0  = 0&& \quad \forall t \in \T_{-1}\\
  & \vI_t^1\vV_t - \vI_{t-1}^1\vV_{t-1} - \vX_{t-1}^1\vV_{t-1}+\ve^{\top}_{t-1} = \v{0} && \quad \forall t \in \T_{-1}\\
    &    s_t^{0+}-I_t^0  +\vh^{\top}\vu_t \ge 0 && \quad \forall t \in \T_{-1}\\
  & \vW^{\top}\vu_t = (\vS_t^{1+}\vV_t-\vI_t^1\vV_t)^{\top} && \quad \forall t \in \T_{-1}\\
  &   s_t^{0-}+I_t^0  + \vh^{\top}\vv_t  \ge 0  && \quad \forall t \in \T_{-1}\\
  & \vW^{\top}\vv_t = (\vS_t^{1-}\vV_t+\vI_t^1\vV_t)^{\top} && \quad \forall t \in \T_{-1}\\
    &   x_t^0  +\vh^{\top}\vgamma_t \ge 0      && \quad \forall t \in \T_{-T}\\
  & \vW^{\top}\vgamma_t \ge  (\vX_t^1\vV_t)^{\top}  && \quad \forall t \in \T_{-T} \\
  &  U^x -x_t^0 + \vh^{\top}\vdelta_t \ge 0  && \quad \forall t \in \T_{-T} \\
  & \vW^{\top}\vdelta_t =  -(\vX_t^1\vV_t)^{\top}  && \quad \forall t \in \T_{-T} \\
  & s_t^{0+}  +\vh^{\top}\vlambda_t \ge 0 && \quad \forall t \in \T_{-1}\\
  & \vW^{\top}\vlambda_t = (\vS_t^{1+}\vV_t)^{\top} && \quad \forall t \in \T_{-1}\\
  &   s_t^{0-}  + \vh^{\top}\vmu_t  \ge 0  && \quad \forall t \in \T_{-1}\\
  & \vW^{\top}\vmu_t = (\vS_t^{1-}\vV_t)^{\top} && \quad \forall t \in \T_{-1}\\
  &  \vu_{t},\ \vv_{t},\ \vlambda_{t},\ \vmu_{t} \in \Re_+^{m}    && \quad \forall t \in \T_{-1} \\
  &   \vgamma_{t},\ \vdelta_{t} \in \Re_+^{m}    && \quad \forall t \in \T_{-T}
  \label{model:inv_SP_multist_LDR_eq_count_end}
	\end{alignat}
\end{subequations} 
where $\vu_{t},\ \vv_{t},\ \vlambda_{t},\ \vmu_{t},\ \vgamma_{t},\ \vdelta_{t}$ are dual variables. We let $ \vgamma_{1} = \vdelta_{1}  = \v{0}$. 

Defining the adaptive decisions in the objective function with the corresponding LDRs, the newsvendor problem's LASC  becomes 
\begin{subequations} \label{model:newsvendor_AASC}
	\begin{alignat}{4}
	  \min ~ &  \sum_{t=1}^{T-1}C_t(x_t^0 +\vX_{t}^1 \vV_t \E[\vd])+ \sum_{t=2}^{T} \left[ H_t (s_t^{0+} +\vS_{t}^{1+} \vV_t \E[\vd]) + B_t(s_t^{0-}+ \vS_{t}^{1-} \vV_t \E[\vd]) \right]   \\
	  \st~~  & \rm{eqs.} \eqref{model:inv_SP_multist_LDR_eq_count_start}-\eqref{model:inv_SP_multist_LDR_eq_count_end}&&
	\end{alignat}
\end{subequations} 
where $\E[\vd]$ is the mean vector of the uncertain demand.

\subsection{Piecewise linear adaptive stochastic counterpart of the newsvendor problem}
In \cite{georghiou2015generalized}, a generic lifting operator is defined as $ L: \Re^k \rightarrow \Re^{k^\prime}$ ($\vd \rightarrow \vd^\prime$), whereas the inverse operator named retraction is  $ R: \Re^{k^\prime} \rightarrow \Re^{k}$ ($\vd^\prime \rightarrow \vd$). The original and lifted uncertainty spaces are defined as $\Xi \subset \Re^k$ and $\Xi^\prime \subset \mathbb{R}^{k^\prime}$($k^\prime > k$), respectively . 

To derive $\vd^\prime$ and $\Xi^\prime$, which is required to define a PLDR, \cite{georghiou2015generalized} first identify the breakpoints where the change of slope occurs. For  $l_t \leq d_t \leq u_t $, we have 
$$ l_t < z_1^t < z_2^t < \dots < z_{r_t -1}^t <u_t \quad t \in \T_{-1}$$
where $z_i^t$ is the $i^{th}$ breakpoint in $d_t$ and $r_t -1$ is the number of breakpoints. The demand in the lifted space $\vd^\prime \in \Re^{k^\prime}$ is given as 
\begin{equation}
	\quad \vd^\prime = (\vd_2^\prime,\dots,\vd_{T}^\prime)^{\top} = (d_{21}^\prime,\dots,d_{2r_2}^\prime,\dots,d_{T1}^\prime,\dots,d_{Tr_{T}}^\prime)^{\top}
\end{equation}
where $k^\prime = \sum_{i=2}^{T} r_t$. Using the set of breakpoints, \cite{georghiou2015generalized} define the lifting operator $L_{ij}(d_i)$ that maps $d_i \in \Xi$ to $d^\prime_{ij} \in \Xi^\prime$ as follows
	\begin{align}
	 L_{ij}(d_i) =
	\begin{cases}
	 d_i \qquad &\text{if}\ r_t = 1,\\
	 \min\{d_i,z_1^t\} \qquad &\text{if}\ r_t > 1,\ j=1,\\
	\max\{\min\{d_i,z_j^i\}-z_{j-1}^i,0\} \qquad &\text{if}\ r_t > 1,\ j=2,\dots,r_t-1\\
	\max\{d_i-z_{j-1}^i,0\} \qquad & \text{if}\ r_t > 1,\ j=r_t,
	\end{cases}\quad \forall i =\{2,\cdots, T\}
	\label{eq:Georghiou_lifting_operator}
	\end{align}
The default case $r_t =1$ corresponds to no lifting. The retraction operator $R_i(\vd^\prime_i)$ is defined as the sum of the lifted elements
\begin{equation}\label{eq:Georghiou_retraction_operator}
d_i = R_i(\vd^\prime_i) = \sum_{j=1}^{r_i} d^\prime_{ij} \quad \forall i =\{2,\cdots, T\}
\end{equation}
The PLDRs, which are LDRs in the lifted uncertainty space, are defined similarly to eq. \eqref{model:inv_SP_multist_LDR_eq}
\begin{subequations} \label{model:inv_SP_multist_Lifted_eq}
	\begin{alignat}{4}
	& x_t (\vd^\prime_{[t]}) =&& ~x_t^{\prime 0} + \vX_t^{\prime 1}\vV_t^\prime\vd^\prime && \qquad \forall  t \in \T_{-T} \\
	& I_t (\vd^\prime_{[t]}) =&& ~I_t^{\prime 0} + \vI_t^{\prime 1}\vV_t^\prime\vd^\prime && \qquad \forall t \in \T \\
	& s_t^+(\vd^\prime_{[t]})=&& ~s_t^{\prime 0+} + \vS_t^{\prime 1+}\vV_t^\prime\vd^\prime && \qquad \forall t \in \T_{-1} \\
	& s_t^-(\vd^\prime_{[t]})=&& ~s_t^{\prime 0-} + \vS_t^{\prime 1-}\vV_t^\prime\vd^\prime && \qquad \forall t \in \T_{-1}
	\end{alignat}
\end{subequations}
where $x_t^{\prime 0},\ I_t^{\prime 0},\ s_t^{\prime 0+},\  s_t^{\prime 0-} \in \Re$ are the intercepts and $\vX_t^{\prime 1},\ \vI_t^{\prime 1},\ \vS_t^{\prime 1+}$, $ \vS_t^{\prime 1-} \in \Re^{1 \times k^\prime}$ are the slopes. We let $x_{1}^{\prime 0}= x_1 $, $I_1^{\prime 0} = I_1$ and $\vX_1^{\prime 1}=\vI_1^{\prime 1}=\textbf{0}$.

The observation matrix $\vV_t^\prime \in \Re^{k^\prime \times k^\prime}$ in the lifted space is reformulated accordingly
\begin{equation}\label{eq:lifted_obs_matrix}
\vd^\prime_{[t]} = \vV_t^\prime \vd^\prime = \left[\begin{array}{ll}
\vI_{k_t } & \mathbf{0}_{k_t \times (k^\prime-k_t) } \\
\mathbf{0}_{(k^\prime-k_t) \times k_t} & \mathbf{0}_{(k^\prime-k_t) \times (k^\prime-k_t)} \\
\end{array}\right] \vd^\prime  = [\vd_2^\prime,\vd_3^\prime,\cdots,\vd_{t}^\prime,0,\dots,0]^{\top} \quad \forall t \in \T_{-1}
\end{equation}
where $k_t=\sum_{i=2}^{t}r_i$ and $k_{T} \equiv k^\prime$. 

\begin{figure}[H]
\centering
\begin{subfigure}{0.26\textwidth}
\centering
\includegraphics[scale=0.4]{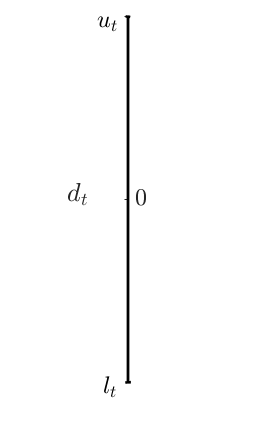}
\caption{No lifting}
\end{subfigure}
\begin{subfigure}{0.36\textwidth}
\centering
\includegraphics[scale=0.4]{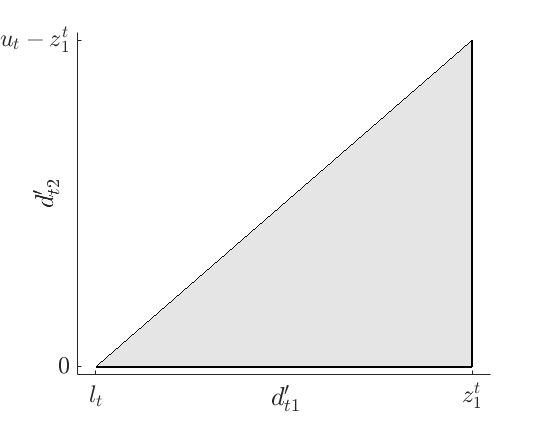} 
\caption{Single breakpoint}
\end{subfigure}
\begin{subfigure}{0.36\textwidth}
\centering
\includegraphics[scale=0.4]{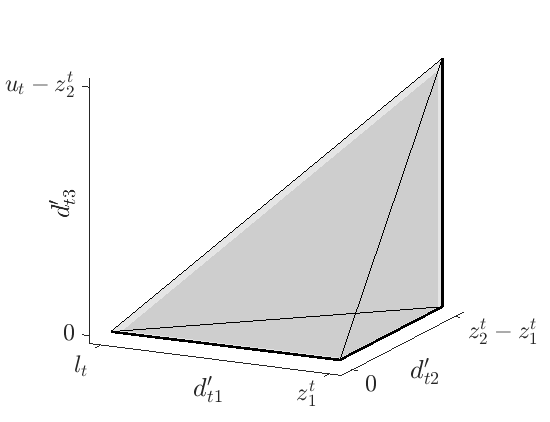} 
\caption{Two breakpoints}
\end{subfigure}
\caption{Lifting the line segment uncertainty set $\Xi_t$  using one and two breakpoints generates a non-convex set $\Xi_t^\prime$ (i.e., bold line). The shaded region represents the corresponding convex hull.}
\label{fig:newsvendor_lifted_conv_hull}
\end{figure}

The lifted demand $\vd^\prime_t$ belongs to a non-convex uncertainty set $\Xi^{\prime}_t$, which violates the strong duality property required to derive the stochastic counterpart. Alternatively, the convex hull of $\Xi^{\prime}_t$ can be used without affecting the optimal solution. It is derived by \cite{georghiou2015generalized} as 
	\begin{align} \label{eq:Georghiou_lifted_convex_hull}
		\begin{split}
		\text{conv}\ \Xi_t^\prime   :&= \{ \vd_t^\prime \in \Re^{r_t}|\  \mathbf{Q}_t^{-1}(1,\vd_t^{\prime \top})^{\top} \geq 0\}\\
		&   =\{ \vd_t^\prime \in \Re^{r_t}|\ \vA_t\vd_t^{\prime} \geq \vb_t\} \qquad \forall t \in \T_{-1}
		\end{split}
	\end{align}
where $\vA_t \in \Re^{(r_t+1)\times (r_t)},\ \vb_t \in \Re^{r_t+1}$ are equal to
$$ \vA_t= \begin{bmatrix}
 -\frac{1}{z_1^t-l_t}&                        & \\
\frac{1}{z_1^t-l_t} & -\frac{1}{z_2^t-z_1^t} &    \\
                  &  \frac{1}{z_2^t-z_1^t} & \ddots \\
 & & \ddots & -\frac{1}{z_{r_t-1}^t-z_{r_t-2}^t}                              \\
  & &        & \frac{1}{z_{r_t-1}^t-z_{r_t-2}^t}  & -\frac{1}{u_t-z_{r_t-1}^t} \\
  & &        &                               &\frac{1}{u_t-z_{r_t-1}^t}      \\
\end{bmatrix},  \quad
\vb_t = \begin{bmatrix}
\frac{z_1^t}{z_1^t-l_t} \\
-\frac{l_t}{z_1^t-l_t}  \\
 0 \\
\vdots \\
\vdots \\
 0\\
\end{bmatrix}
$$
Figure \ref{fig:newsvendor_lifted_conv_hull} illustrates the original uncertainty set $\Xi_t$, the non-convex lifted uncertainty set $\Xi^\prime_t$ and its convex hull using one and two breakpoints, respectively.

The dimension of the lifted uncertain parameter in stage $t$ is analogous to the number of scenarios in stage $t$ in a scenario tree. Figure~\ref{fig:analogy_scenario_decision_rule} illustrates the aforementioned analogy which supports the motivation behind relating an empirical efficient construct of a scenario tree to the lifting strategy in a hybrid decision rule.

The convex hull of a lifted generic polytope $\Xi$ as in eq.~\eqref{eq:multistg_inv_polyhedral_uncertainty_set} has no tractable representation. \cite{georghiou2015generalized} present a tractable outer approximation of the convex hull which is the intersection of (i) the convex hull of the lifted uncertainty set of the smallest hyper-rectangle containing the generic polytope as in eq.~\eqref{eq:Georghiou_lifted_convex_hull} and (ii) the convex lifted uncertainty set of the generic polytope $\Xi$. The outer approximation is given as
\begin{align} \label{eq:polyhedral_lifted_convex_hull}
		\begin{split}
\bar{\Xi}^\prime :& = \{\vd^\prime \in \Re^{k^\prime}|\ \v{W}\ R\vd^\prime \geq \v{h},\  \vA_t\vd_t^\prime \geq \vb_t \quad t \in \T_{-1} \}  \\
           & = \{\vd^\prime \in \Re^{k^\prime}|\ \vW^\prime \vd^\prime \geq \v{h},\  \vA_t\vd_t^\prime \geq \vb_t \quad t \in \T_{-1} \} \\
           & =  \{ \vd^\prime \in \Re^{k^\prime}|\ \vA^{\rm{l}} \vd^\prime \geq  \vb^{\rm{l}} \}  \quad \supseteq \text{conv }\Xi^\prime          
	\end{split}
\end{align}
The matrix $\vW^\prime \in \Re^{T-1 \times k^\prime}$ is defined as
$$
\vW^\prime= \begin{bmatrix} \vw_{1}\v1_{r_2}, \vw_{2}\v1_{r_3},\cdots,\vw_{T-1}\v1_{r_{T}} 
\end{bmatrix}
$$
where $\vw_i \in \Re^{T-1}$ is the $i^{th}$ column of $\vW$, and $\v{1}_{r_i} \in \R^{1\times r_1}$ is a unit row vector. The matrix $\vA^{\rm{l}} \in \Re^{(m+m^\prime)\times k^\prime}$ and vector $\vb^{\rm{l}} \in \Re^{m+m^\prime}$, where $m^\prime = \sum_{t=2}^{T}(r_t+1)$, are formulated as

$$ \vA^{\rm{l}}= \begin{bmatrix}
 & & \vW^\prime&                  &    \\
 \vA_2   & \v{0}     & \cdots & \v{0}       \\
  \v{0}      & \vA_3 & \cdots     & \v{0}       \\
  \vdots &  \vdots & \ddots   &   \vdots              \\
  \v{0}      &   \v{0}    &    \cdots & \vA_{T}\\
\end{bmatrix},  \quad
\vb^{\rm{l}} = \begin{bmatrix}
\vh \\
\vb_2  \\
\vb_3 \\
\vdots \\
\vb_{T}\\
\end{bmatrix}
$$

Since both uncertainty sets $\Xi$ and $\bar{\Xi}^\prime$ are polytopes, the procedure of constructing the Piecewise linear adaptive stochastic counterpart (PWLASC) is similar to that of LASC with few changes in parameters and dimensions. The formulation of the multistage newsvendor's PWLASC is given as 
\begin{subequations} \label{model:inv_SP_multist_Lifted_finite_polyhedral}
	\begin{alignat}{4}\nonumber
	 \min~~ & \sum_{t=1}^{T-1}C_t(x_t^{\prime 0} +\vX_{t}^{\prime 1} \vV^\prime_t \E[\vd^\prime])+ \sum_{t=2}^{T}(H_t (s_t^{\prime 0+}+\vS_{t}^{\prime 1+} \vV^\prime_t \E[\vd^\prime]) + B_t(s_t^{\prime 0-}+ \vS_{t}^{\prime 1-} \vV^\prime_t \E[\vd^\prime]) )   \\ 
	 \st~~  & \rm{eqs.}\ \eqref{model:inv_SP_multist_LDR_eq_count_start}-\eqref{model:inv_SP_multist_LDR_eq_count_end}   
	\end{alignat}
\end{subequations}
where $\vW = \vA^{\rm{l}},\ \vh = \vb^{\rm{l}},\ \vV_t = \vV_t^\prime,\ \ve_{t-1} = \ve^\prime_{t-1},\ m  \rightarrow m +m^\prime$, and $\E[\vd^\prime]$ is the mean vector of the lifted uncertain demand. 
The row vector $\ve^\prime_{t} \in \Re^{k^\prime}$ has a value of 1 from index $\sum_1^{t-1}r_t +1$ to $\sum_1^{t}r_t$ and a value of 0 elsewhere.

\subsection{Numerical results for the multistage stochastic newsvendor problem}
\label{sec:newsvendor_numerical_results}
This section illustrates solution quality improvements generated by PLDRs, corroborating the empirical results found in  \cite{georghiou2015generalized}. In the next section, we address the following question: Does the additional flexibility of PLDRs (i.e., their higher resolution) always lead to better practical policies?  We demonstrate that an improvement is not guaranteed; in fact, the solution can deteriorate. Unless otherwise stated, the independent uncertain demand in each stage follows a uniform distribution between 0 and 10. The cost coefficients are fixed in all stages at $ C_t = 3,\ H_t = 1.5,\  B_t = 7$. The ordering amount limit $U^x$ and the initial inventory $I_1$ are equal to 8 and 4, respectively.

Three decision rules are studied: (i) LDR, which serves as a lower bound, (ii) PLDR-1 ($\E[d_t]$) and (iii) PLDR-1 ($U^x$). The last two PLDRs apply one breakpoint in $d_t$ in each stage at the mean value and at the ordering amount limit, respectively. While there is not a systematic method to identify the optimal set of breakpoints, we think that relating it to a characteristic of the uncertainty distribution ($\E[d_t]$) or to a physical parameter of the system ($U^x$) is an intuitive and a practical choice.

For $T$ = 4, the optimal cost generated by an LDR is equivalent to 83.5, and the first stage solution $x_1$ is equal to 8. As expected, a PLDR-1 ($\E[d_t]$) reduces the optimal cost to 66.25 which reflects a 20.66$\%$ decrease from the LDR case. Likewise, $x_1$ has decreased by 25$\%$ to a value of 6. Meanwhile, a PLDR-1 ($U^x$) solution quality exceeds that of a PLDR-1 ($\E[d_t]$) where the optimal cost is equal to 63.60, and $x_1$ is reduced to 4. The optimal policies for the three decision rules are shown in Table \ref{tab:optimal_policies_multistage_inventory_solution_T_4} in Appendix A.
 
\begin{figure}[H]
\begin{subfigure}{0.32\textwidth}
\centering
\includegraphics[width=5cm,height=4cm]{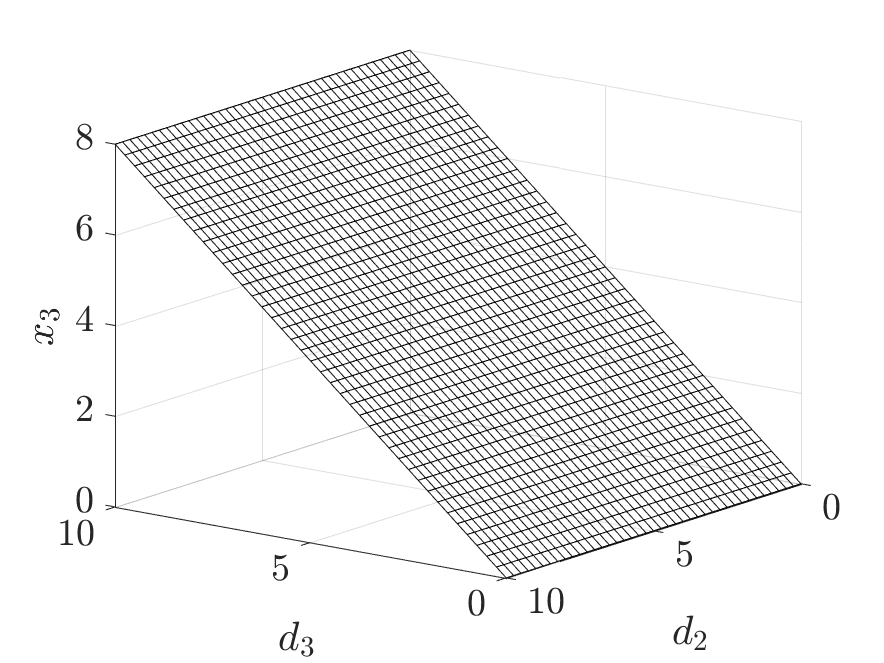}
\end{subfigure}
\begin{subfigure}{0.32\textwidth}
\centering
\includegraphics[width=5cm,height=4cm]{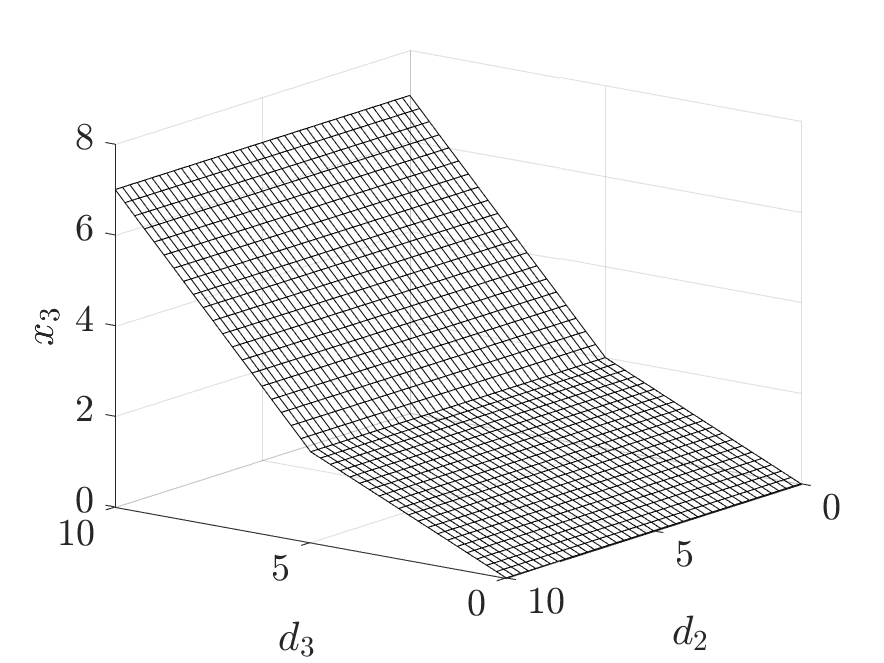}
\end{subfigure}
\begin{subfigure}{0.32\textwidth}
\centering
\includegraphics[width=5cm,height=4cm]{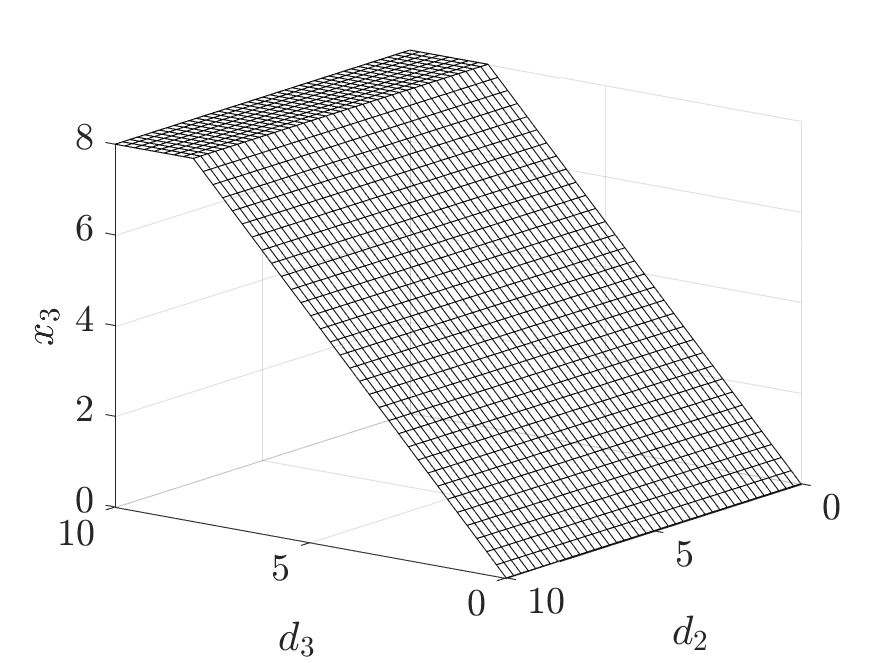}
\end{subfigure} 
\\
\begin{subfigure}{0.32\textwidth}
\centering
\includegraphics[width=5cm,height=4cm]{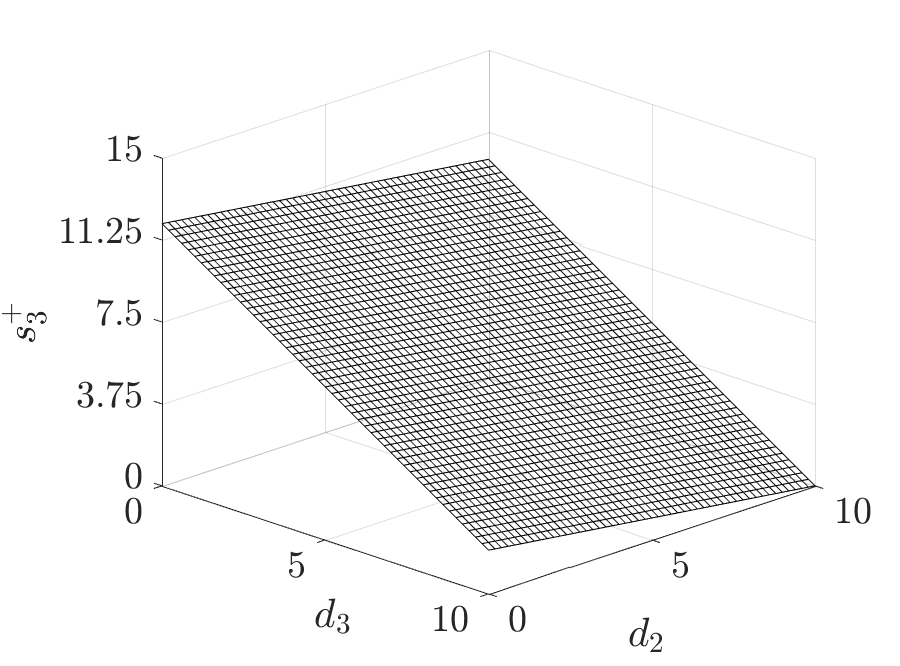}
\end{subfigure}
\begin{subfigure}{0.32\textwidth}
\centering
\includegraphics[width=5cm,height=4cm]{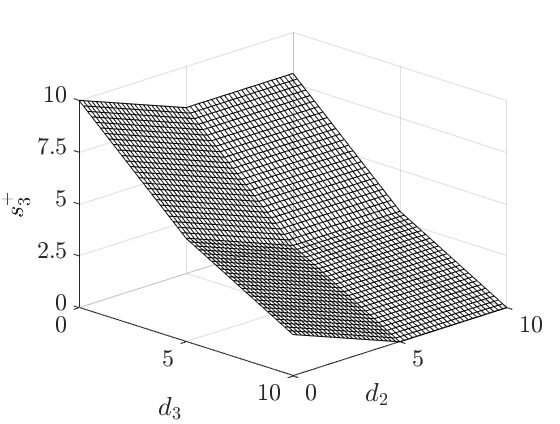}
\end{subfigure}
\begin{subfigure}{0.32\textwidth}
\centering
\includegraphics[width=5cm,height=4cm]{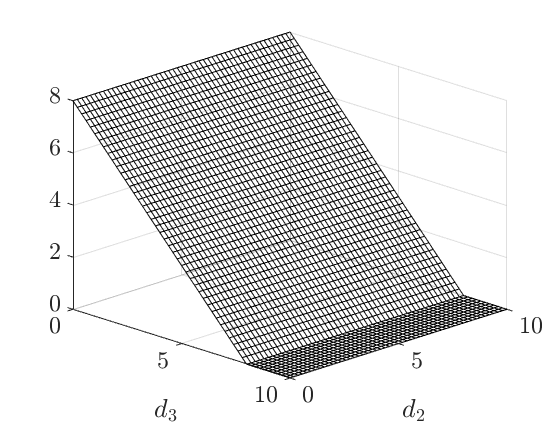}
\end{subfigure}
\\
\begin{subfigure}{0.32\textwidth}
\centering
\includegraphics[width=5cm,height=4cm]{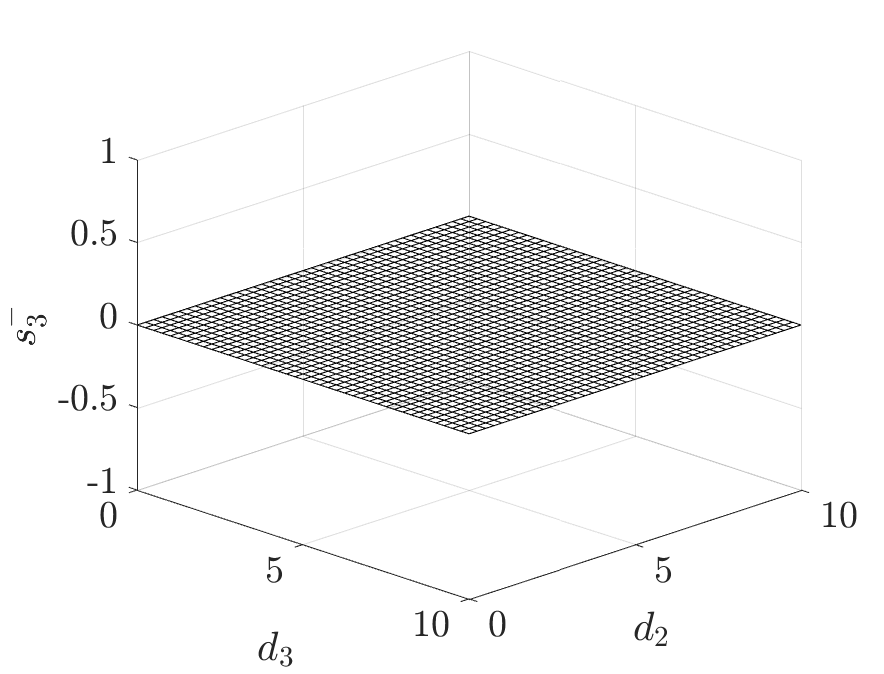}
\caption{}
\end{subfigure}
\begin{subfigure}{0.32\textwidth}
\centering
\includegraphics[width=5cm,height=4cm]{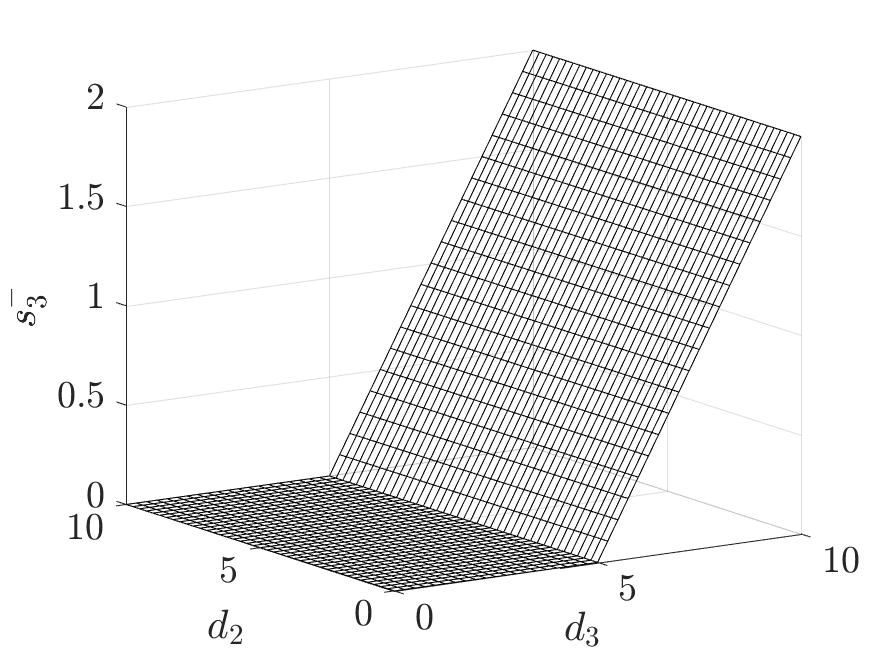}
\caption{}
\end{subfigure}
\begin{subfigure}{0.32\textwidth}
\centering
\includegraphics[width=5cm,height=4cm]{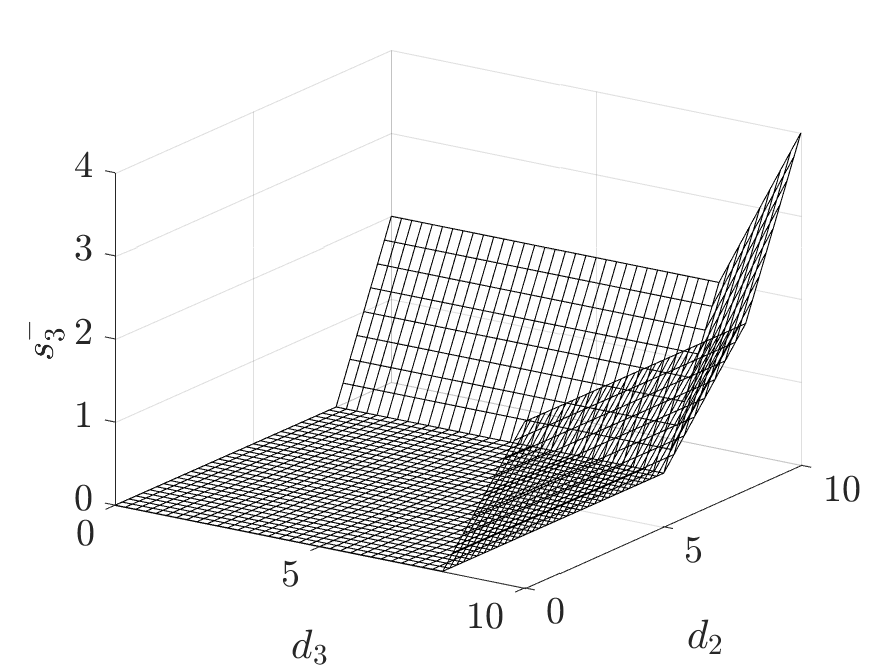}
\caption{}
\end{subfigure}
\caption{Profiles of $x_3(d_2,d_3),\ s_3^+(d_2,d_3)\ \text{and}\  s_3^-(d_2,d_3)$ using (a) an LDR, (b) a PLDR-1 ($\E[d_t]$) and (c) a PLDR-1 ($U^x$). The presence of a single breakpoint in $d_1$ and $d_2$ generates at most four recourse strategies as in Figures (b) and (c). The recourse rate of the ordering policy for PLDR-1 ($U^x$) is the highest at a value of 1, reaching the ordering limit $U^x$ at $d_2=8$. Further, note that the optimal policies in Figures (b) and (c) output overestimated decisions for specific regions in the $(d_2,d_3)$ space. The states variables $s_3^+$ and $s_3^-$ simultaneously have positive values when this should never be the case in any applied policy since by definition they are mutually exclusive. Parameters used: $T=4,\ U^x=8,\ I_1=4,\ d_2,d_3 \sim \mathcal{U}(0,10)$.  }
\label{fig:inv_illustrattion_pol}
\end{figure}

The ordering, inventory and backlog policies in stage 3 are demonstrated in Figure \ref{fig:inv_illustrattion_pol}. PLDRs provide the flexibility to implement different recourse strategies based on the previous realization of the uncertain demand. For example, $s_3^+(d_2,d_3)$ policy generated by a PLDR-1 ($\E[d_t]$) exhibits four different recourse strategies due to a single breakpoint in $d_2$ and $d_3$, while an LDR generates only a single recourse strategy.

\section{Can an LDR ever be superior to a PLDR?} 
In some cases, LDR policies are proven to be optimal as in \cite{bertsimas2010Iancu}, otherwise one would expect a PLDR to outperform an LDR because it is more flexible.  In this section, we show that this need not be the case even when both decision rules are constructed from the same underlying uncertainty set. This example bolsters the argument of simulating policies (i.e., decision rules) within a simulator as opposed to merely assessing their objective function values in a look-ahead model. In our study, the look-ahead models are the derived adaptive stochastic counterparts. For brevity, we will refer to the ``look-ahead model'' as the ``model'' for the rest of the paper. For a thorough discussion on the motivation to assess the quality of decision policies obtained in a look-ahead model via a simulator, see \cite{powell2014clearing}.

In this section and throughout the study, we compute all ``optimal" decision rules at the outset of the problem. Then, we gradually evaluate the adaptive decisions as subsets of the uncertainty are sequentially revealed. This is also known as \textit{closed loop} policy/implementation and is the basis of our comparison. We do not study decision rules implemented in an \textit{open loop} or rolling horizon manner where new static and adaptive decisions are recursively computed  for the entire planning horizon after implementing the previously obtained optimal static decisions, observing the first revealed subset of the uncertainty and shifting the problem one stage forward.

First, we present an example of overestimated inventory and backlog decisions observed in Figure \ref{fig:inv_illustrattion_pol}. For ($d_2,d_3$)=($7,9$), the values of ($s_3^+(d_2,d_3)$,$s_3^-(d_2,d_3)$) are equivalent to $(0.6,1.6)$ and $(1,1)$ using a PLDR-1 ($\E[d_t]$) and a PLDR-1 ($U^x$), respectively. This consequently leads to an overestimated model-based cost.

Next, we assess the solution quality via a simulator using $10^5$ samples. We implement only the optimal ordering policy and compute the state variables by their definitions: $s_t^+ = \max(0,I_t)$ and $s_t^- = \max(0,-I_t)\ \forall t \in \T_{-1}$. A comparison between the model- and simulator-based costs, for the three policies in Figure \ref{fig:inv_illustrattion_pol}, is demonstrated in Table \ref{tab:cost_comp_solver_simulator}. The discrepancy between the model- and simulator-based costs alters the superiority of a DR with respect to another DR. For example, a PLDR-1 ($U^x$) appears superior to a PLDR-1 ($\E[d_t]$) using the model, however this is not the case when evaluated within a simulator. Consequently, the set of breakpoints that is ought to be optimal in a model, may not be optimal within a simulator environment.

\begin{table}[H]
  \centering
  \caption{Comparison between the model- and simulator-based costs for the multistage stochastic newsvendor problem using an LDR, a PLDR-1 ($\E[d_t]$) and a PLDR-1 ($U^x$). The model-based cost exhibits a degree of overestimation which may be misleading in terms of the quality of a decision rule as shown for PLDR-1 ($U^x$) and PLDR-1 ($\E[d_t]$). This advocates assessing the quality of the optimal policies within a simulator.} \label{tab:cost_comp_solver_simulator}
    \begin{tabular}{lccccc}
    \toprule
       & Model-based & \multicolumn{4}{c}{Simulator-based}\\
       	\cmidrule(lr){3-6}
     Decision rule     & $\E[\rm{cost}]$ & $\E[\rm{cost}]$ & $\sigma$ & min & max    \\
    \midrule 
    LDR                  & 83.50 & 75.14 & 4.72 & 63.00    & 88.20    \\
    PLDR-1 ($\E[d_t]$)  & 66.25 & 59.07 & 9.97 & 45.30 & 104.40   \\
    PLDR-1 ($U^x$)      & 63.60 & 59.88 & 11.23& 36.52 & 139.44  \\
    \bottomrule
    \end{tabular}
  \end{table}
The reason for this outcome is due to the presence of mutually exclusive state variables and the  robust nature of the stochastic counterparts. For example, the stochastic counterparts ensure feasibility for the worst uncertainty realization, consequently they guarantee feasibility for any other possible uncertainty realization. However, enforcing complementary slackness of inventory and backlog decision rules (i.e., state variables) for the worst uncertainty realization \textit{does not} guarantee complimentary slackness at any other possible uncertainty realizations. As a result, the look-ahead model-based policies may output simultaneous positive decisions values for what are supposed to be in practice mutually exclusive decisions (i.e, the seller either has a deficit or surplus of the good). 

\begin{figure}[H]
\centering
\begin{subfigure}{0.48\textwidth}
\centering
\includegraphics[width=7.5cm,height=5cm]{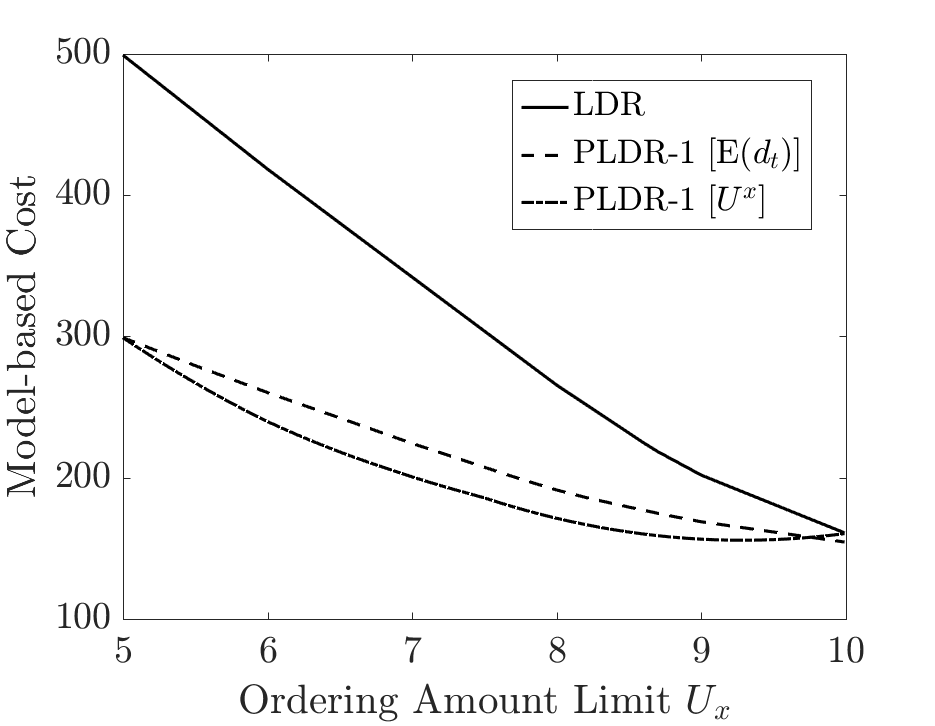}
\end{subfigure}
\begin{subfigure}{0.48\textwidth}
\centering
\includegraphics[width=7.5cm,height=5cm]{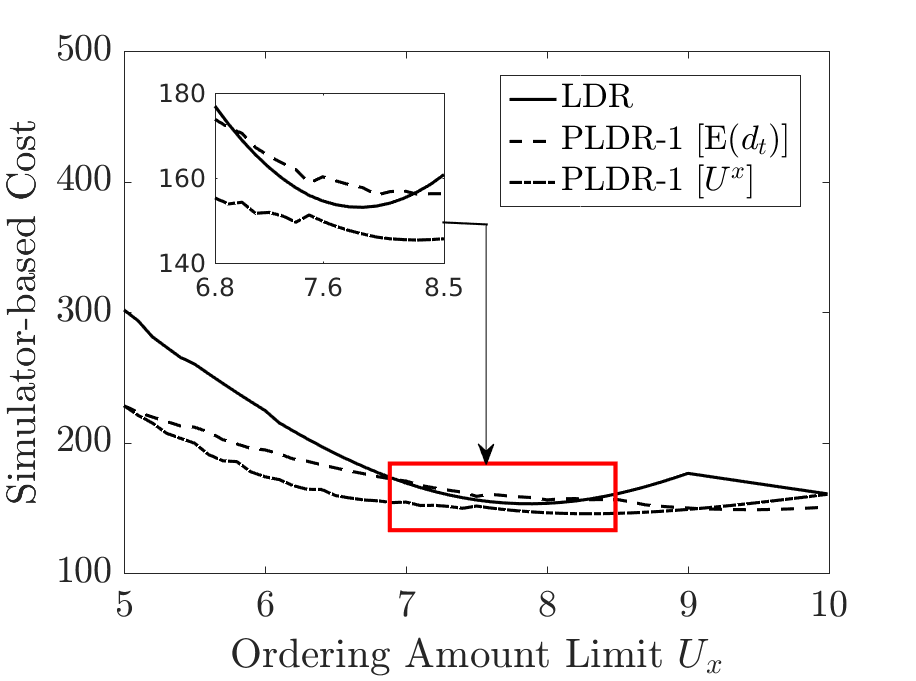}
\end{subfigure}
\caption{ Profiles of the model- and simulator-based costs using an LDR, a PLDR-1 ($\E[d_t]$) and a PLDR-1 ($U^x$) as a function of $U^x$ for the multistage stochastic newsvendor problem. As expected, PLDRs-1 are always superior to LDRs with respect to the model-based cost. However, for $U^x \in [\text{{\raise.17ex\hbox{$\scriptstyle\sim$}}}6.9,\text{{\raise.17ex\hbox{$\scriptstyle\sim$}}}8.3]$, the simulator-based cost generated by the LDR policy is counter-intuitively superior to that of PLDR-1 ($\E[d_t]$). Parameters used: $T=8,\ I_1=4,\ d_t \sim \mathcal{U}(0,10)\ t\in\T_{-1}$.}
\label{fig:newsvendor_comp_solver_simulator}
\end{figure}

\begin{figure}[H]
\centering
\begin{subfigure}{0.48\textwidth}
\centering
\includegraphics[width=7.5cm,height=5cm]{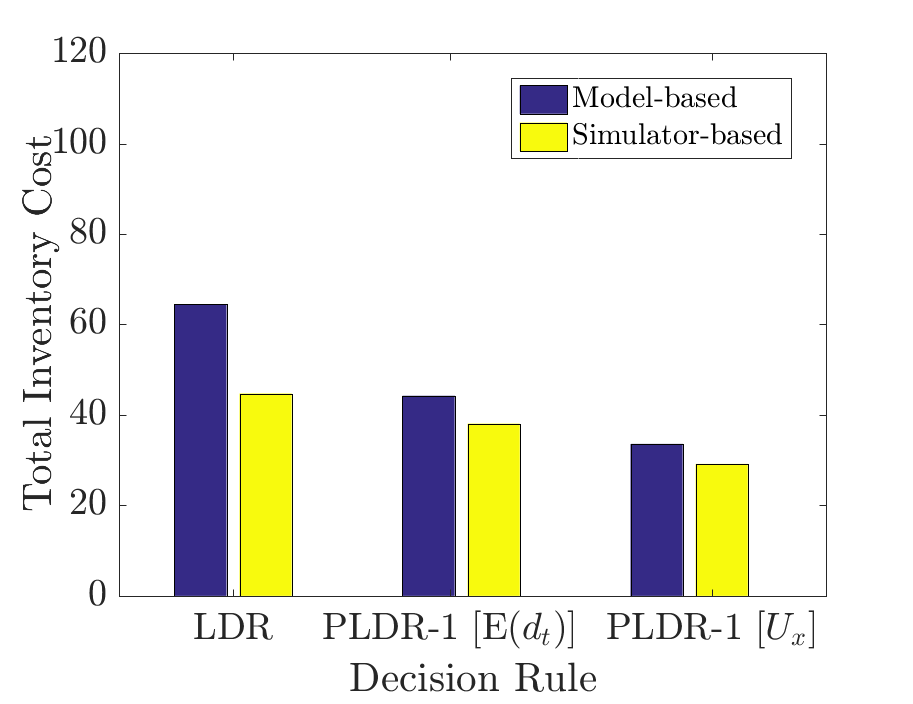}
\end{subfigure}
\begin{subfigure}{0.48\textwidth}
\centering
\includegraphics[width=7.5cm,height=5cm]{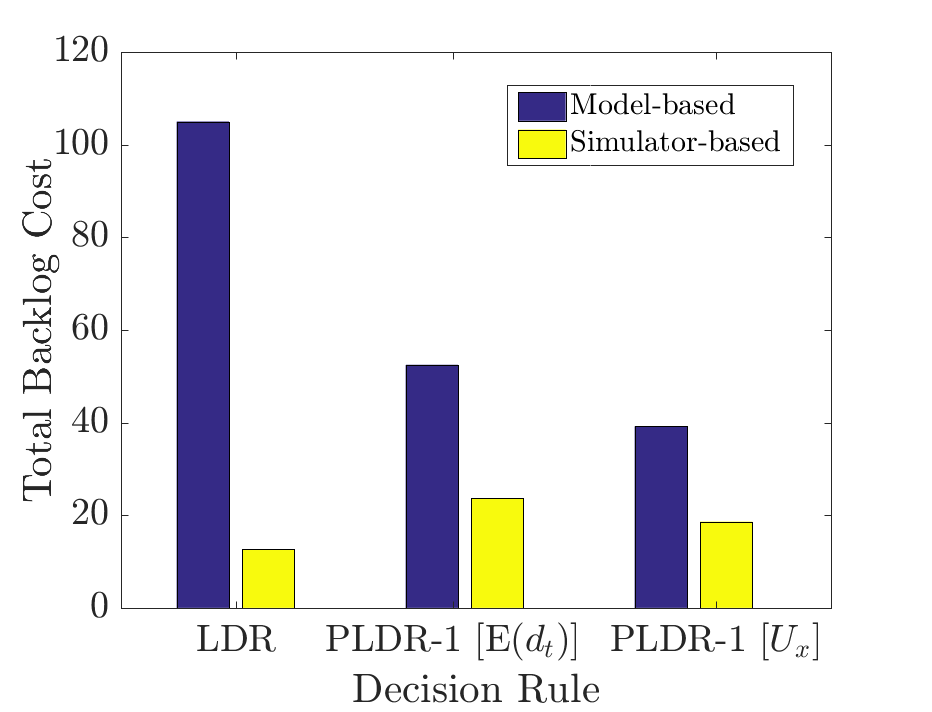}
\end{subfigure}
\caption {Total inventory and backlog costs of a multistage stochastic newsvendor problem using an LDR, a PLDR-1 ($\E[d_t]$) and a PLDR-1 ($U^x$) via the model and simulator. The overestimation is predominant in the backlog cost component, and it is most prominent in the LDR policy. In contrast,  simulated backlog cost for the LDR is less than that for both PLDRs. This is due to the conservative LDR ordering policy which explains why LDR outperforms PLDR-1($\E[d_t]$) at this specific computational setting. Parameters used: $T=8,\ U^x=8,\  I_1=4,\ d_t \sim \mathcal{U}(0,10)\ t\in\T_{-1}$. } \label{fig:news_com_str_bkl_solv_sim}
\end{figure}

Now, the question is: With respect to the simulator-based cost, can an LDR policy outperform a PLDR policy? Extending the planning horizon $T$ from 4 to 8, the model- and simulator-based costs are computed over $U^x \in [5,10]$ for an LDR, a PLDR-1($\E[d_t]$) and a PLDR1 ($U^x$); the profiles are shown in Figure~\ref{fig:newsvendor_comp_solver_simulator}.  Indeed, a PLDR-1 outperforms an LDR in terms of the model-based cost regardless of the single breakpoint's value. However, with respect to the simulator-based cost, an LDR is found to be superior to a PLDR-1 ($\E[d_t]$) for $U^x \in [\text{{\raise.17ex\hbox{$\scriptstyle\sim$}}}6.9,\text{{\raise.17ex\hbox{$\scriptstyle\sim$}}}8.3]$. Statistically, using $10^5$ samples and $100$ replications, the magnitude of the simulated cost variance for the three decision rules is found to be in the order of $1e^{-26}$ at a fixed value of $U_x$ (e.g., 7). Further, the same trend emerges in the simulator-based cost for the three decision rules while using $10^4$, $10^3$ and $10^2$ samples.   
To justify this observation and appreciate the overestimation in the model-based cost, Figure \ref{fig:news_com_str_bkl_solv_sim} compares the total inventory and backlog costs for both the model and simulator at $T=8$ and $U^x= 8$. Note that since the ordering policy is implementable, the model- and simulator-based total ordering costs are equal for any DR at any computational setting. For the DRs investigated, the overestimation is most significant in the backlog cost component, in particular it is more prominent in the LDR policy. Despite the fact that the model-based backlog cost exhibited by the LDR is the highest, the corresponding simulator-based cost is the lowest in comparison to the two PLDRs investigated. This is due to the more conservative ordering policy, and it explains why the LDR is superior to the PLDR-1 ($\E[d_t]$) within the simulator in this specific computational setting.
As a final observation, it is worth to note that the uncertainty set used in Figure \ref{fig:newsvendor_comp_solver_simulator} is a hyper-rectangle. The convex hull of the aforementioned lifted uncertainty set is constructed exactly with no overestimation/approximation (see eq.~\eqref{eq:polyhedral_lifted_convex_hull}). Thus the deterioration in the flexibility of a PLDR-1($\E[d_t]$) in a simulator environment, consequently its inferiority to an LDR, can not be attributed to an overestimation in the lifted uncertainty set. 

\section{Hybrid lifting strategies}
In this section we exploit the modularity of PLDRs with axial segmentation. Using a multistage stochastic transportation problem, we empirically illustrate that having higher uncertainty resolution (i.e., more linear pieces) in early stages is more significant than having it in late stages. Then, we demonstrate that non-increasing HDRs are more competitive than non-decreasing HDRs in terms of the trade-off between solution quality and computational time. But, as in any design problem, we show that a poorly designed non-increasing HDR loses the computational benefits. 

\subsection{Multistage stochastic transportation model}
fWe now present a multistage stochastic transportation problem. The framework involves a set of $\I$ suppliers and a set of $\J$ customers. Suppliers are able to store some of their products for future time as inventory.   In each stage $t$, a decision maker must determine (1) the number of units $x_{it}$ produced by each supplier $i$ which will first become available for delivery in the subsequent stage $t+1$; (2) the number of transported units $y_{ijt}$ from supplier $i$ to costumer $j$; and (3) the number of units carried as inventory $I_{it}$ by each supplier $i$.

The formulation of the deterministic multistage transportation problem is
\begin{subequations} \label{model:multistage_deterministic_transportation_problem}
\begin{alignat}{4} \label{eq:obj_trnsprt}
\max_{x_{it},I_{it},y_{ijt}}~~& \sum_{t \in \T_{-1}} \sum_{i \in \I} \sum_{j \in \J} (R_{jt} - T_{ijt}) y_{ijt}
- \sum_{t \in \T_{-T}} \sum_{i \in \I} C_{it} x_{it}
-  && \sum_{t \in \T_{-1}} \sum_{i \in \I} H_{it} I_{it} + \sum_{i \in \I} S_i I_{iT} && \\
\st ~~& \sum_{j \in \J} y_{ijt} \leq I_{it} && \quad \forall i \in \I,\ t \in \T_{-1} \label{eq:supply_leq_inv} \\
& \sum_{i \in \I} y_{ijt} \leq D_{jt} && \quad \forall j \in \J,\ t \in \T_{-1} \label{eq:supply_leq_demand} \\
& I_{it} = I_{it-1} + x_{it-1} - \sum_{j \in \J} y_{ijt-1} && \quad \forall i \in \I,\ t \in \T_{-1} \label{eq:inv_balance_deterministic} \\
\label{eq:trnspt_uxi_deterministic}
& 0 \le x_{it} \leq U_i^{\max} && \quad \forall i \in \I,\ t \in \T_{-T} \\
& I_{it},\ y_{ijt} \geq 0 & & \quad \forall i \in \I,\ j \in \J,\ t \in \T_{-1} 
\end{alignat}
\end{subequations}
where $y_{ij1} = 0\ \forall i,j$ and $I_{i1}$ is the initial inventory.

The objective function \eqref{eq:obj_trnsprt} seeks to maximizes the profit. It includes revenue, transportation, production, holding costs and salvage value to mitigate the end of horizon effect \citep{fisher2001ending}. Constraint \eqref{eq:supply_leq_inv} ensures that the amount transported out from each supplier does not exceed the amount of inventory currently available. Meanwhile, constraint \eqref{eq:supply_leq_demand} dictates that the amount transported to each customer does not surpass the demand. The inventory balance is satisfied by eq. \eqref{eq:inv_balance_deterministic} and equation \eqref{eq:trnspt_uxi_deterministic} sets the production limit for each supplier $i$.

Equation \eqref{eq:supply_leq_inv} can be defined in terms of the initial inventory $I_{i1}$ and the cumulative difference between produced and transported amounts by supplier $i$ 
\begin{equation} \label{eq:inv_balance_reform}
I_{it} = I_{i1} + \sum_{k=1}^{t-1} (x_{ik} - \sum_{j \in \J} y_{ijk}) \quad \forall i \in \I, t \in \T_{-1}
\end{equation}
Model \eqref{model:multistage_deterministic_transportation_problem} can be simplified by substituting $I_{it}$ in eq. \eqref{eq:supply_leq_inv} by the right hand side of eq.~\eqref{eq:inv_balance_reform}, and dropping out eq. \eqref{eq:inv_balance_deterministic}. For our study, we use Model \eqref{model:multistage_deterministic_transportation_problem} without any modifications.

The uncertainty is represented by a primitive vector $\vxi = \{\xi_2,\cdots,\xi_{T}\}$, which the uncertain demand $D_{jt}$ is a function of. Introducing $\vxi$ into Model \eqref{model:multistage_deterministic_transportation_problem},  the multistage stochastic transportation problem is given as
\begin{subequations} \label{model:multistage_stochastic_transportation_problem_generic_formulation_2}
\begin{alignat}{4}
\max_{\substack{x_{it}(\cdot),I_{it}(\cdot)\\ y_{ijt}(\cdot)}}~~& \E \Bigg[ \sum_{t \in \T_{-1}} \sum_{i \in \I} \sum_{j \in \J} (R_{jt} - T_{ijt}) y_{ijt}(\vxi_{[t]})- \sum_{t \in \T_{-T}}&& \sum_{i \in \I}  C_{it} x_{it}(\vxi_{[t]})  \notag \\
	  &   - \sum_{t \in \T_{-1}} \sum_{i \in \I} H_{it}I_{it}(\vxi_{[t]}) + \sum_{i \in \I} S_i I_{iT}(\vxi_{[T]}) \Bigg] && \\
\st ~~& \sum_{j \in \J} y_{ijt}(\vxi_{[t]}) \leq I_{it}(\vxi_{[t]}) &&  \forall  i \in \I,\ t \in \T_{-1},\ \vxi \in \Xi  \\
& \sum_{i \in \I} y_{ijt}(\vxi_{[t]}) \leq D_{jt}(\xi_{t}) && \forall  j \in \J,\ t \in \T_{-1},\ \vxi \in \Xi \\
& I_{it}(\vxi_{[t]}) = I_{it-1}(\vxi_{[t-1]}) + x_{it-1}(\vxi_{[t-1]}) -  \sum_{j \in \J} &&y_{ijt-1}(\vxi_{[t-1]})\notag \\
&  && \forall  i \in \I,\ t \in \T_{-1},\ \vxi \in \Xi \\\label{eq:trnsp_xit_upp}
& 0 \le x_{it}(\vxi_{[t]}) \leq U_i^{\max} && \forall  i \in \I,\ t \in \T_{-T},\ \vxi \in \Xi \\
& I_{it}(\vxi_{[t]}),\  y_{ijt}(\vxi_{[t]}) \geq 0 && \forall  i \in \I,\ j \in \J,\ t \in \T_{-1},\ \vxi \in \Xi
\end{alignat}
\end{subequations}
where $\Xi$ is the underlying polyhedral uncertainty set (see eq.~\eqref{eq:multistg_inv_polyhedral_uncertainty_set}). We let $y_{ij1}(\vxi_{[1]}) = y_{ij1}$, $x_{i1}(\vxi_{[1]}) = x_{i1}$ and $I_{i1}(\vxi_{[1]}) = I_{i1}$. 
\subsubsection{Linear adaptive stochastic counterpart of transportation problem}
LDRs of the adaptive decisions are defined as
\begin{subequations} \label{model:trns_SP_multist_LDR_eq}
	\begin{alignat}{4}
	& x_{it} (\vxi_{[t]})   =&& ~x_{it}^0 + \vX_{it}^1\vV_t\vxi && \qquad \forall i \in \I,\  t \in \T_{-T} \\
	& y_{ijt} (\vxi_{[t]})  =&& ~y_{ijt}^0 + \vY_{ijt}^1\vV_t\vxi && \qquad \forall i \in \I,\ j \in \J,\ t \in \T \\
    & I_{it} (\vxi_{[t]})   =&& ~I_{it}^0 + \vI_{it}^1\vV_t\vxi && \qquad \forall i \in \I,\  t \in \T
	\end{alignat}
\end{subequations}
where $x_{it}^0,\ y_{ijt}^0,\ I_{it}^0 $ are the intercepts and $\vX_{it}^1,\ \vY_{ijt}^{1},\ \vI_{it}^1 \in \Re^{1 \times (T-1)}$ are the slopes. We let $x_{i1}^0 = x_1,\ y_{ij1}^0 = 0,\ I_{i1}^0 = I_{i1}\ \text{and}\ \vY_{ij1}^1 = \vX_{i1}^1 =\vI_{i1}^1=\v{0}$.

The customers' uncertain demand in stage $t$ is assumed to be a linear function of $\xi_t$
\begin{equation}
D_{jt}(\xi_t) = D_{jt}^0 + D_{jt}^1\xi_t \qquad \forall j \in \J,\ t\in \T_{-1}
\end{equation}
where $D_{jt}^0\ \text{and}\ D_{jt}^1$ are parameters. Implementing LDRs in Model \eqref{model:multistage_stochastic_transportation_problem_generic_formulation_2}, we obtain
\begin{subequations} \label{model:multistage_stochastic_transportation_problem_LDR}
\begin{alignat}{4}
\max_{\substack{x^0_{it},\ \vX^1_{it},\ I^0_{it}\\ \vI^1_{it},\ y^0_{ijt},\ \vY^1_{ijt}}}~~& \E \Bigg[ \sum_{t \in \T_{-1}} \sum_{i \in \I} \sum_{j \in \J} (R_{jt} - T_{ijt}) (y_{ijt}^0 + \vY_{ijt}^1\vV_t\vxi) && \notag \\
 &- \sum_{t \in \T_{-T}} \sum_{i \in \I}  C_{it} (x_{it}^0 + \vX_{it}^1\vV_t\vxi) - \sum_{t \in \T_{-1}} \sum_{i \in \I} H_{it}&& (I_{it}^0 + \vI_{it}^1\vV_t\vxi) \notag \\
 & + \sum_{i \in \I} S_i (I_{iT}^0 + \vI_{iT}^1\vV_T\vxi) \Bigg] && \\
		\st ~~& \sum_{j \in \J} (y_{ijt}^0 + \vY_{ijt}^1\vV_t\vxi) \leq I_{it}^0 + \vI_{it}^1\vV_t\vxi && \forall  i \in \I,\ t \in \T_{-1},\ \vxi \in \Xi  \\
		& \sum_{i \in \I} (y_{ijt}^0 + \vY_{ijt}^1\vV_t\vxi) \leq D_{jt}^0 + D_{jt}^1\xi_{t} && \forall  j \in \J,\ t \in \T_{-1},\ \vxi \in \Xi \\
		& I_{it}^0 + \vI_{it}^1\vV_t\vxi =I_{it-1}^0 + \vI_{it-1}^1\vV_{t-1}\vxi + x_{it-1}^0 +&& \vX_{it-1}^1 \vV_{t-1}\vxi  \notag \\
		&- \sum_{j \in \J} (y_{ijt-1}^0 + \vY_{ijt-1}^1\vV_{t-1}\vxi) &&  \forall  i \in \I,\ t \in \T_{-1},\ \vxi \in \Xi \\
		& x_{it}^0 + \vX_{it}^1\vV_t\vxi \leq U_i^{\max} && \forall  i \in \I,\ t \in \T_{-T},\ \vxi \in \Xi \\
		& x_{it}^0 + \vX_{it}^1\vV_t\vxi \geq 0 && \forall  i \in \I,\ t \in \T_{-T},\ \vxi \in \Xi \\
		& I_{it}^0 + \vI_{it}^1\vV_t\vxi \geq 0 &&  \forall  i \in \I,\ t \in \T_{-1},\ \vxi \in \Xi \\
		& y_{ijt}^0 + \vY_{ijt}^1\vV_t\vxi \geq 0 && \forall  i \in \I,\ j \in \J,\ t \in \T_{-1},\ \vxi \in \Xi
\end{alignat}
\end{subequations}
Similar to the procedure followed for the newsvendor problem, the overall transportation problem's LASC is given as

\begin{subequations} \label{model:transportation_affine_counterpart} 
	\begin{alignat}{4}
	 \max_{\substack{x^0_{it},\ \vX^1_{it},\ I^0_{it}\\ \vI^1_{it},\ y^0_{ijt},\ \vY^1_{ijt}}}~~&   \sum_{t \in \T_{-1}} \sum_{i \in \I} \sum_{j \in \J} (R_{jt} - T_{ijt}) (y_{ijt}^0 + \vY_{ijt}^1\vV_t\E[\vxi]) -  && \sum_{t \in \T_{-T}} \sum_{i \in \I}  C_{it} (x_{it}^0 + \vX_{it}^1\vV_t\E[\vxi])  \notag \\
 & - \sum_{t \in \T_{-1}} \sum_{i \in \I} H_{it}(I_{it}^0 + \vI_{it}^1\vV_t\E[\vxi]) + \sum_{i \in \I} S_i (I_{iT}^0 +  && \vI_{iT}^1\vV_T\E[\vxi]) \notag  \\ 
 \label{model:trns_SP_multist_LDR_eq_count_start}
	   &   I_{it}^0 -\sum_{j \in \J} y_{ijt}^0 +\vh^{\top}\vu_{it} \ge 0 && \quad \forall i \in \I,\ t \in \T_{-1}\\    
  & \vW^{\top}\vu_{it} = (\vI_{it}^1\vV_t-\sum_{j \in \J} \vY_{ijt}^1\vV_t)^{\top} && \quad \forall i \in \I,\ t \in \T_{-1}\\   
  &   D_{jt}^0-\sum_{i \in \I} y_{ijt}^0 + \vh^{\top}\vv_{jt}  \ge 0  &&\quad  \forall j \in \J,\ t \in \T_{-1}\\
  & \vW^{\top}\vv_{jt} = ( D_{jt}^1\ve_{t-1}-\sum_{j \in \J} \vY_{ijt}^1\vV_t)^{\top} && \quad \forall j \in \J,\ t \in \T_{-1}\\
  & I_{it}^0 - I_{it-1}^0 - x_{it-1}^0 + \sum_{j \in \J} y_{ijt-1}^0 = 0 && \quad \forall i \in \I,\ t \in \T_{-1}\\
  & \vI_{it}^1\vV_t - (\vI_{it-1}^1+\vX_{it-1}^1-\sum_{j \in \J} \vY_{ijt-1}^1)\vV_{t-1}  =\textbf{0} && \quad \forall i \in \I,\ t \in \T_{-1} \\
   & U_i^{\max} -x_{it}^0 +\vh^{\top}\vmu_{it} \ge 0 &&\quad \forall i \in \I,\ t \in \T_{-T}\\
  & \vW^{\top}\vmu_{it} = -(\vX_{it}^{1}\vV_t)^{\top} &&\quad \forall i \in \I,\ t \in \T_{-T}\\
  & x_{it}^0 +\vh^{\top}\vlambda_{it} \ge 0 &&\quad \forall i \in \I,\ t \in \T_{-T}\\
  & \vW^{\top}\vlambda_{it} = (\vX_{it}^{1}\vV_t)^{\top} &&\quad \forall i \in \I,\ t \in \T_{-T}\\
  &  I_{it}^0 + \vh^{\top}\vdelta_{it} \ge 0  && \quad \forall i \in \I,\ t \in \T_{-1} \\
  & \vW^{\top}\vdelta_{it} =  (\vI_{it}^1\vV_t)^{\top}  && \quad \forall i \in \I,\ t \in \T_{-1} \\
   & y_{ijt}^0 +\vh^{\top}\vgamma_{ijt} \ge0 &&\quad \forall i \in \I,\ j \in \J,\ t \in \T_{-1}\\
  & \vW^{\top}\vgamma_{ijt} = (\vY_{ijt}^{1}\vV_t)^{\top} &&\quad \forall  i \in \I,\ j \in \J,\ t \in \T_{-1}\\
  & \vlambda_{it}, \vmu_{it} , \vu_{it}, \vv_{jt}, \vdelta_{it}, \vgamma_{ijt} \in \Re_+^{m}  && \quad \forall i \in \I,\ j \in \J,\ t \in \T_{-1} 
  \label{model:trns_SP_multist_LDR_eq_count_end}
	\end{alignat}
\end{subequations} 
where $\E[\vxi]$ is the mean vector with respect to the distribution of $\vxi$ and $\vu_{it},\ \vv_{jt},\ \vlambda_{it},\ \vmu_{it},\ \vdelta_{it},\ \vgamma_{ijt}$ are dual variables. We let $\vlambda_{i1}= \vmu_{i1}=\v{0}$.

\subsubsection{Piecewise Linear adaptive stochastic counterpart of the transportation problem}

The adaptive recourse decisions are defined as PLDRs (i.e., LDRs in the lifted space)
\begin{subequations} \label{model:trns_SP_multist_Lifted_eq}
	\begin{alignat}{4}
	& x_{it} (\vxi^\prime_{[t]})   =&& ~x_{it}^{\prime 0} + \vX_{it}^{\prime 1}\vV^\prime_t\vxi^\prime && \qquad \forall i \in \I,\  t \in \T_{-T} \\
	& y_{ijt} (\vxi^\prime_{[t]})  =&& ~y_{ijt}^{\prime 0} + \vY_{ijt}^{\prime 1}\vV^{\prime}_t\vxi^\prime && \qquad \forall i \in \I,\ j \in \J,\ t \in \T \\
    & I_{it} (\vxi^\prime_{[t]})   =&& ~I_{it}^{\prime 0} + \vI_{it}^{\prime 1}\vV^{\prime}_t\vxi^\prime && \qquad \forall i \in \I,\  t \in \T
	\end{alignat}
\end{subequations}
where $x_{it}^{\prime 0},\ y_{ijt}^{\prime 0},\ I_{it}^{\prime 0} $ are the intercepts and $\vX_{it}^{\prime 1},\ \vY_{ijt}^{\prime 1},\ \vI_{it}^{\prime 1} \in \Re^{1 \times k^\prime}$ are the slopes. We let $x_{i1}^{\prime 0} = x_1,\ y_{ij1}^{\prime 0} = 0,\ I_{i1}^{\prime 0} =I_{i1}\ \text{and}\ \vX_{i1}^{\prime 1} = \vY_{ij1}^{\prime 1} = \vI_{i1}^{\prime 1}=\v{0}$. 

Using the retraction operator defined in eq. \eqref{eq:Georghiou_retraction_operator}, the customers' uncertain demand is reformulated as a function of $\vxi_t^\prime = \{\xi^\prime_{11},\cdots,\xi^\prime_{1r_t}\}$
\begin{equation}\label{eq:trnps_affine_demand}
D_{jt}(\vxi_t^\prime) = D_{jt}^0 + D_{jt}^1\sum_{i=1}^{r_t}\xi^\prime_{it}, \qquad \forall j \in \J,\ t\in \T_{-1}
\end{equation}

The transportation problem's PWLASC is similar to the LASC with changes in parameters and dimensions
\begin{subequations} \label{model:trans_LASC}
	\begin{alignat}{4}
	  \max_{\substack{x^{\prime 0}_{it},\ \vX^{\prime 1}_{it},\ I^{\prime 0}_{it}\\ \vI^{\prime 1}_{it},\ y^{\prime 0}_{ijt},\ \vY^{\prime 1}_{ijt}}}~~&   \sum_{t \in \T_{-1}} \sum_{i \in \I} \sum_{j \in \J} (R_{jt} - T_{ijt}) (y_{ijt}^{\prime 0} + \vY_{ijt}^{\prime 1}\vV^\prime_t\E[\vxi^\prime]) - \sum_{t \in \T_{-T}} \sum_{i \in \I}  C_{it} (x_{it}^{\prime 0} + \vX_{it}^{\prime 1}\vV^\prime_t\E[\vxi^\prime])&& \notag \\
 & - \sum_{t \in \T_{-1}} \sum_{i \in \I} H_{it}(I_{it}^{\prime 0} + \vI_{it}^{\prime 1}\vV^\prime_t\E[\vxi^\prime]) + \sum_{i \in \I} S_i (I_{iT}^{\prime 0} + \vI_{iT}^{\prime 1}\vV^\prime_T\E[\vxi^\prime])  && \\
	  \st~~  & \rm{eqs.} \eqref{model:trns_SP_multist_LDR_eq_count_start}-\eqref{model:trns_SP_multist_LDR_eq_count_end}&& 
	\end{alignat}
\end{subequations} 
where $\vW = \vA^{\rm{l}},\ \vh = \vb^{\rm{l}},\ \vV_t = \vV_t^\prime,\ \ve_{t-1} = \ve^\prime_{t-1},\ m \rightarrow m +m^\prime $, and $\E[\vxi^\prime]$ is the mean vector with respect to the distribution of $\vxi^\prime$. The parameters $\vA^{\rm{l}}$ and  $ \vb^{\rm{l}}$ define the outer approximation of the convex hull of the lifted uncertainty set as in eq.~\eqref{eq:polyhedral_lifted_convex_hull}. The observation matrix in the lifted space $\vV_t^\prime$ is given in eq.~\eqref{eq:lifted_obs_matrix}.

\subsection{Numerical results for the multistage stochastic transportation problem}
\label{sec:transporttion_numerical_results}
In this section, we demonstrate that HDRs with higher uncertainty resolution in early stages are more flexible than HDRs with higher uncertainty resolution in late stages. This complements with the empirical evidence found in scenario-based stochastic programming methods where scenario trees with higher granularity in early stages are more attractive (\cite{bakkehaug2014stochastic}, \cite{arslan2017bulk}). We first investigate the sensitivity of the solution quality with respect to uncertainty resolution in each stage, then we design a set of hybrid decision rules to demonstrate the acquired computational benefits.
 
Unless otherwise stated, our computational setting includes three suppliers $\I= \{1,2,3\}$ and two customers $\J= \{1,2\}$. The uncertainty in each stage is independent and follows a uniform distribution between 0 and 3. The remaining parameters are listed in Table \ref{tab:param-trnsprt}.
\begin{table}[H]
\caption{Parameters defining the computational setting of a multistage stochastic transportation problem. They are assumed to be constant for all stages (e.g. $C_{it} \equiv C_{i}\ \forall t$).}
\begin{center}
\begin{tabular}{lccccccccccccc}
\toprule
  & $C_{i}$ & $H_{i}$  &  $S_{i} $&$U_i^{\max}$& & &\multicolumn{2}{c}{$T_{ij}$}&  & &$D^0_{j}$ & $D^1_{j}$  & $R_{j}$  \\
\cmidrule(lr){8-9}
$i \downarrow$    & & & & & & $i \downarrow j \rightarrow$&  1 & 2 & & $j \downarrow $    & & &   \\
\cmidrule(lr){1-5}  \cmidrule(lr){7-9} \cmidrule(lr){11-14}
 1   & 5 & 2  & 0 & 10& & 1  &3 & 4  & &1   & 5 & 3 & 18\\
 2   & 7 & 3  & 0 & 8 & & 2 &1 & 5   & & 2  & 2 & 1 & 16 \\
 3   & 1 & 0.5& 0 & 5 &   & 3& 6 & 2  & & ---& --- & --- & ---\\
\bottomrule
\end{tabular}
\label{tab:param-trnsprt}
\end{center}

\end{table}

\subsubsection*{Limitation of PLDRs with single lifting component}
\begin{figure}[H]
\begin{subfigure}{0.48\textwidth}
\centering
\includegraphics[width=7.5cm,height=5cm]{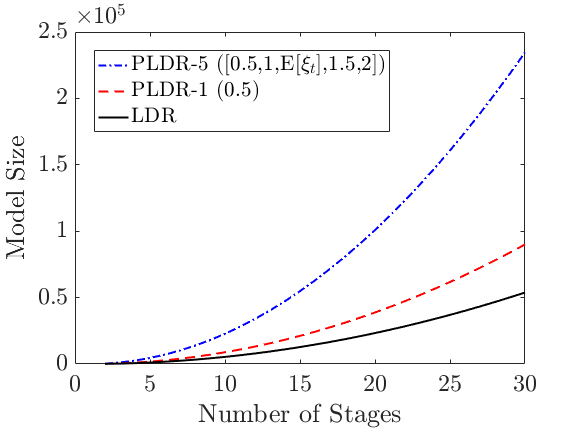}
\end{subfigure}
\begin{subfigure}{0.48\textwidth}
\centering
\includegraphics[width=7.5cm,height=5cm]{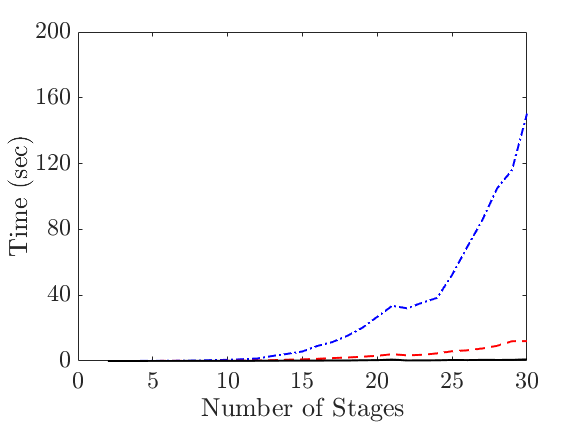}
\end{subfigure}
\caption{Model size (number of constraints and variables) and computational time required to solve a multistage stochastic transportation problem via an LDR, a PLDR-1 (0.5) and a PLDR-5 ([0.5,1,$\E[\xi_t]$,2,2.5]) for $T=\{2,\dots,30\}$. The exponential increase in computational time induced by PLDR-5 limits its applicability in large scale problems and motivates the need for hybrid decision rules that combine the salient features of both LDRs and PLDRs. Parameters used: $ \xi_t \sim \mathcal{U}(0,3)\ \forall t \in \T_{-1}$.  }
\label{fig:tranprt_SP_sens_T_exp_inc_time}
\end{figure}

The improvement in solution quality induced by implementing PLDRs may be significant. Nonetheless, the computational overhead resulting from lifting the uncertainty to a high resolution impedes its applicability in large scale problems. For a range of planning horizons, the model size (number constraint and variables) of the transportation problem, and the computational time using an LDR, a PLDR-1 (0.5) and a PLDR-5 ([0.5,1,$\E[\xi_t]$,2,2.5]) are illustrated in Figure \ref{fig:tranprt_SP_sens_T_exp_inc_time}. The exponential increase in computational time for a PLDR-5 is evident. 

\subsubsection*{Impact of uncertainty resolution in early and late stages on solution quality}
In this section, we perform a sensitivity analysis to investigate the impact of uncertainty resolution (i.e., number of linear pieces) in each stage on the solution quality. To serve as reference points, we first define two base decision rules: (1) LDR where there is no lifting of $\xi_t$ in any stage, and (2) PLDR-5 ([0.5,1,$\E[\xi_t]$,2,2.5]) which is assumed to approximate the true solution. The two base DRs have the least and highest uncertainty resolution, respectively. The sensitivity of solution quality to resolution of $\xi_t$ is computed by varying the number of breakpoints implemented in $\xi_t$ while keeping the number of breakpoints used to lift $\xi_{t^\prime}$ for $ t^{\prime} \in \T_{-\{1,t\}} $  equal to that of the base DR. 

\begin{figure}[H]
\begin{subfigure}{0.96\textwidth}
\includegraphics[width=7.5cm,height=5cm]{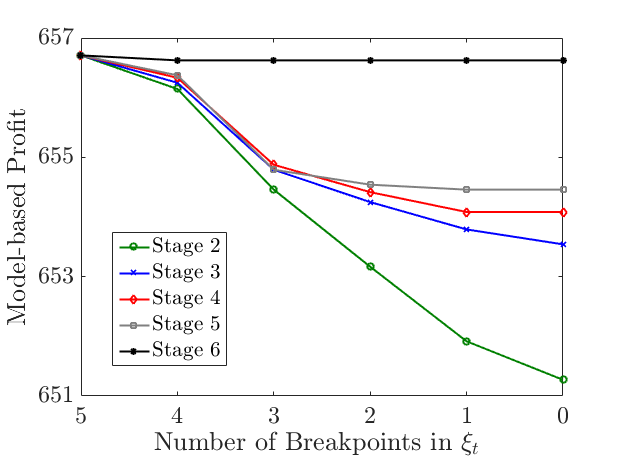}
\includegraphics[width=7.5cm,height=5cm]{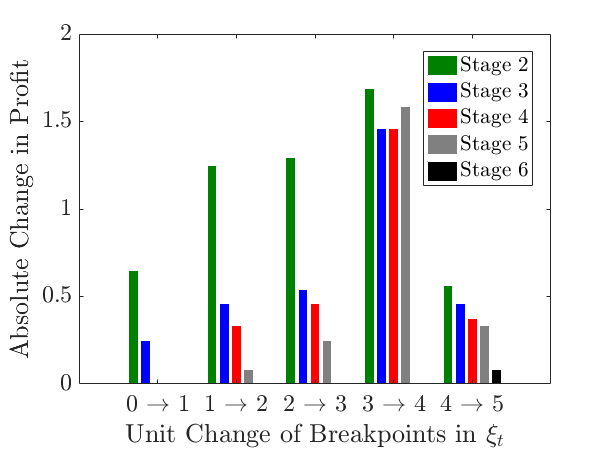}
\caption{Base decision rule: PLDR-5 ([0.5,1,$\E[\xi_t]$,2,2.5]). }
\end{subfigure}
\begin{subfigure}{0.96\textwidth}
\includegraphics[width=7.5cm,height=5cm]{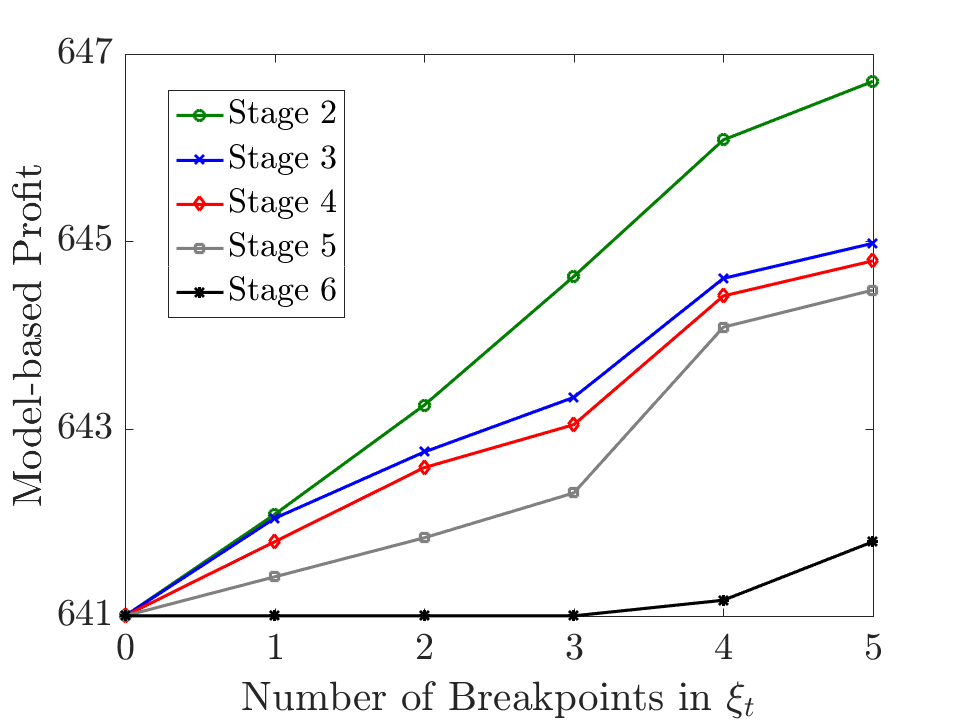}
\includegraphics[width=7.5cm,height=5cm]{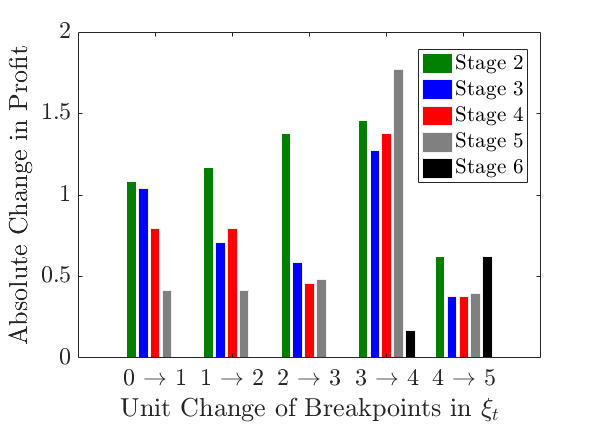}
\caption{Base decision rule: LDR.}\label{fig:secondT6Val6}
\end{subfigure}
\caption{Sensitivity of the model-based profit with respect to the uncertainty resolution in a multistage stochastic transportation problem using a PLDR-5 ($\mathcal{Z}_{\textrm{base}}$) and an LDR as base DRs. Number of breakpoints $i$ corresponds to lifting $\xi_t$ with the first $i$ elements in $\mathcal{Z}_{\textrm{base}}$, while keeping the resolution of $\xi_{t^\prime}$ for $t^\prime \in  \T_{-\{1,t\}}$ at that of the base DR. Higher uncertainty resolution in early stages has a higher impact on the solution quality than in late stages. A similar trend is seen in the bar plots despite some discrepancies as the salvage value is computed empirically (i.e., end of horizon effect). Parameters used: $T =6,\ \mathcal{Z}_{\textrm{base}}=\{0.5,1,\E[\xi_t],2,2.5\},\ S_i=6\ \forall i,\ \xi_t \sim \mathcal{U}(0,3)\ \forall t \in \T_{-1}$. } \label{fig:sens_curves_trns_T6}
\end{figure}

\begin{figure}[H]
\begin{subfigure}{0.96\textwidth}
\includegraphics[width=7.5cm,height=5cm]{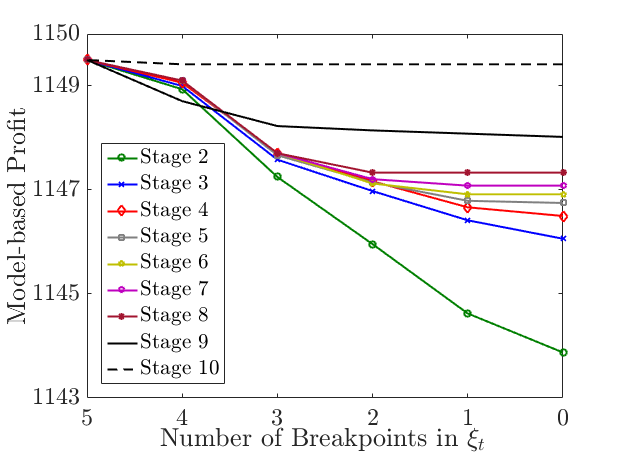}
\includegraphics[width=7.5cm,height=5cm]{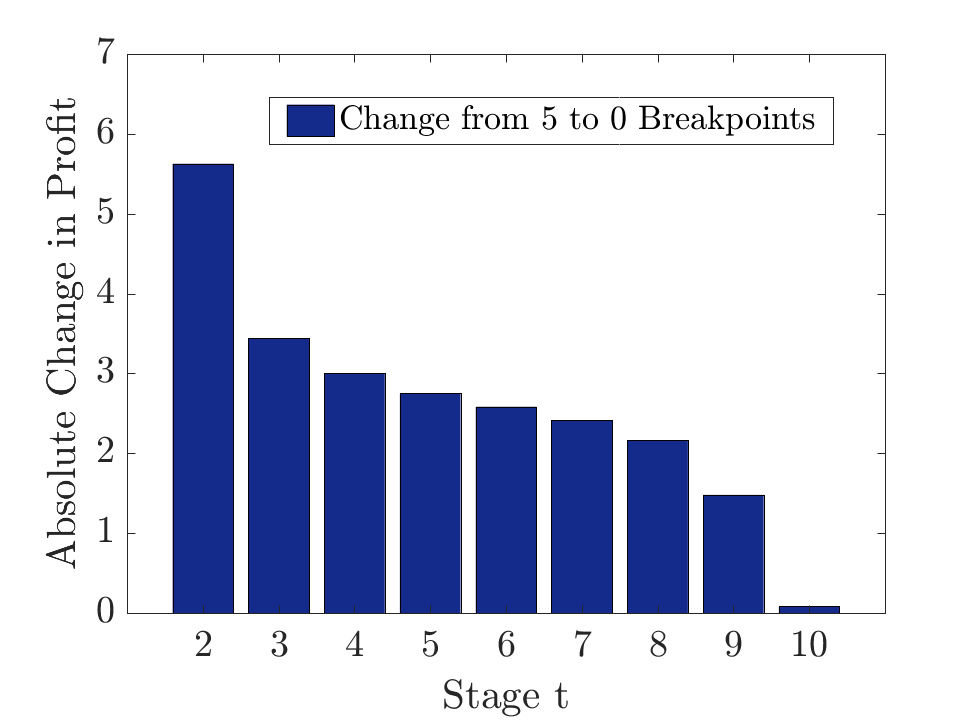}
\caption{Base decision rule: PLDR-5 ([0.5,1,$\E[\xi_t]$,2,2.5]). }
\end{subfigure}

\begin{subfigure}{0.96\textwidth}
\includegraphics[width=7.5cm,height=5cm]{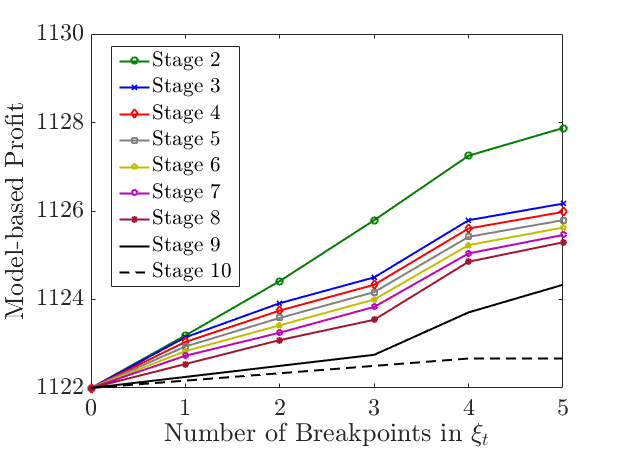}
\includegraphics[width=7.5cm,height=5cm]{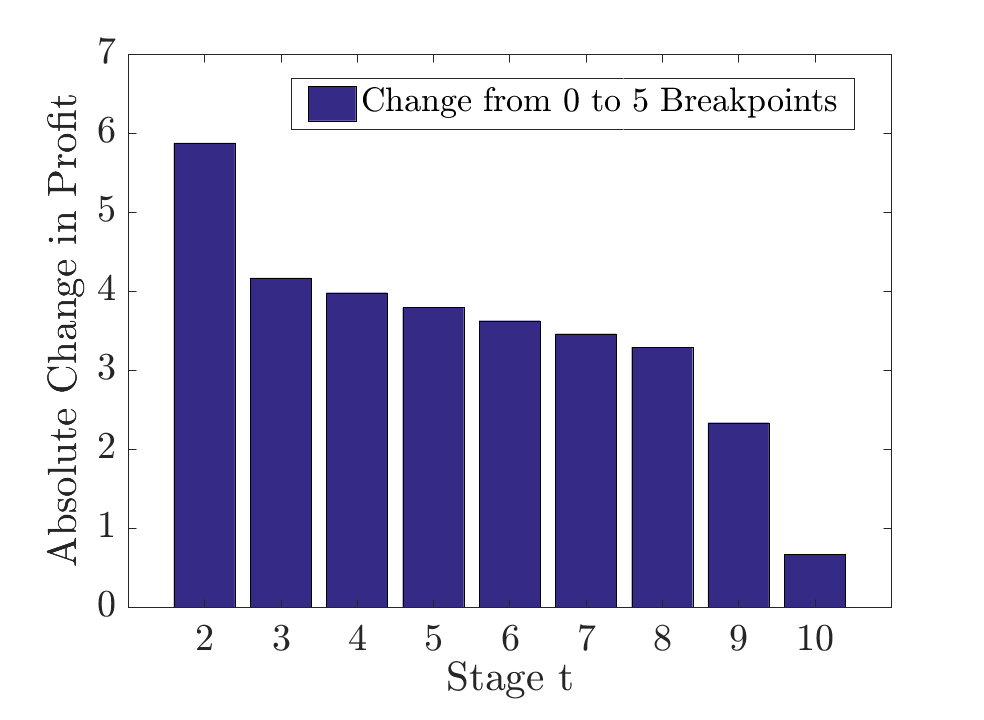}
\caption{Base decision rule: LDR.}
\end{subfigure}

\caption{Sensitivity of the model-based profit with respect to the uncertainty resolution in the multistage stochastic transportation problem using a PLDR-5 ($\mathcal{Z}_{\textrm{base}}$) and an LDR as base DRs. Number of breakpoints $i$ corresponds to lifting $\xi_t$ with the first $i$ elements in $\mathcal{Z}_{\textrm{base}}$, while keeping the resolution of $\xi_{t^\prime}$ for $t^\prime \in  \T_{-\{1,t\}}$ at that of the base DR. Higher uncertainty resolution in early stages has a higher impact on the solution quality than in late stages. This is more evident in the absolute change in profit due to the change of breakpoints between the two extrema 0 and 5. Parameters used: $T=10,\ \mathcal{Z}_{\textrm{base}}=\{0.5,1,\E[\xi_t],2,2.5\},\ S_i=7.5\ \forall i,\ \xi_t \sim \mathcal{U}(0,3)\ \forall t \in \T_{-1}$.  }
\label{fig:sens_curves_trns_T10}
\end{figure}

Setting $T=6$ and $S_i =6$ for all $i$, the sensitivity curves for $\xi_t$ are shown in Figure \ref{fig:sens_curves_trns_T6} for the two base DRs. An adjacent bar plot illustrates the absolute change in profit with respect to the unit change in number of breakpoints as well. It is observed that the sensitivity of the solution quality (i.e., profit) to the change in uncertainty resolution in earlier stages is higher than that in late stages. Even though this is not seen at some instances as in the bar plots in Figure~\ref{fig:sens_curves_trns_T6}(b), we should keep in mind that (1) the salvage value is determined empirically and does not entirely eliminate the end of horizon effect and (2) the impact of lifting $\xi_t$ using an additional breakpoint is highly dictated by the breakpoint value and the nature of the true solution of adaptive decision rules in stage $t$.
Next, we extend the planning horizon $T$ to $10$ and adjust the salvage value to $S_i=7.5$ for all $i$. The sensitivity curves for $\xi_t$ are generated in Figure \ref{fig:sens_curves_trns_T10}. Implementing an LDR as a base DR, the expected trend is clearly obtained where having higher uncertainty resolution in early stages offers higher solution quality. Likewise, for PLDR-5($[0.5,1,\E[\xi_t],2,2.5]$) as a base DR, the deterioration in the solution quality is  highest when the uncertainty resolution is reduced in early stages. There is one discrepancy where reducing the number of breakpoints to $4$ in the $9^{th}$ stage has a higher impact than the previous stages. However, the absolute change in profit induced by changing the number of breakpoints between the two extrema ($0$ and $5$) for the two base DRs clearly indicates that higher resolution in early stages is more attractive, in particular in stage 2 (i.e., in the first revealed instance of the uncertain parameter).

\subsubsection*{Implementing HDRs for a 10-stage planning horizon }
Higher uncertainty resolution in early stages offers more flexibility for PLDR policies compared to the same resolution in late stages. Still, we need to explore various lifting strategies (which is itself a design problem), and answer the following question: Is it guaranteed that non-increasing HDRs will provide the most attractive trade-off between solution quality and computational time? 

To address this question, a set of experiments is conducted to a modified version of the multistage stochastic transportation problem. First, an expansion decision is added where the decision maker decides in the first stage whether or not to invest in additional capacity for each supplier $i$. This decision is represented by the binary variable $y^{\textrm{bin}}_{i}=\{0,1\}$. Second, the produced amount within the additional capacity $x^{\rm{exp}}_{it}$ incurs a $50\%$ increase in the per unit production cost. The objective function is reformulated as follows
\begin{align}
\max_{\substack{x_{it}(\cdot),x^{\rm{exp}}_{it}(\cdot),I_{it}(\cdot)\\ y_{ijt}(\cdot),\ y^{\textrm{bin}}_{i}}}~~& \E \Bigg[ \sum_{t \in \T_{-1}} \sum_{i \in \I} \sum_{j \in \J} (R_{jt} - T_{ijt}) y_{ijt}(\vxi_{[t]}) - \sum_{t \in \T_{-T}} \sum_{i \in \I}  C_{it}(x_{it}(\vxi_{[t]})+(1+\epsilon) x^{\rm{exp}}_{it}(\vxi_{[t]}) )\notag \\
	  &  -  \sum_{i \in \I}  M_{i}y^{\textrm{bin}}_{i} - \sum_{t \in \T_{-1}} \sum_{i \in \I} H_{it}I_{it}(\vxi_{[t]}) + \sum_{i \in \I} S_i I_{iT}(\vxi_{[T]}) \Bigg]  
\end{align} 
where $M_i$ is the capital cost required for the expansion in supplier $i$ and $\epsilon=0.5$ . The upper bound constraint of $x^{\rm{exp}}_{it}\in \R_+$ is given as
\begin{equation}
x^{\rm{exp}}_{it}(\vxi_{[t]}) \le y^{\textrm{bin}}Q_i \quad \forall  \vxi \in \Xi,\ i \in \I,\ t \in \T_{-T}
\end{equation}
where $Q_i$ is the maximum added production capacity for supplier $i$. The inventory balance for all $ i \in \I,\ t \in \T_{-1},\ \text{and}\ \vxi \in \Xi$ is reformulated to take into account both types of production decisions
\begin{equation}
I_{it}(\vxi_{[t]}) = I_{it-1}(\vxi_{[t-1]}) + x_{it-1}(\vxi_{[t-1]}) + x^{\rm{exp}}_{it-1}(\vxi_{[t-1]})-  \sum_{j \in \J} y_{ijt-1}(\vxi_{[t-1]}) 
\end{equation}

Further, the uncertainty follows a uniform distribution between 0 and 1 in each stage: $\xi_t \sim \mathcal{U}(0,1)$. The number of suppliers and customers are increased to 10, and the remaining parameters defining the computational setting are tabulated in Appendix B. The potential breakpoints used to define PLDRs belong to a ``base set'' which is defined as function of the mean ($\mu$) and standard deviation ($\sigma$) of $\mathcal{U}(0,1)$: $\mathcal{Z}_{\textrm{base}} =\{\mu-\sigma,\mu-0.5\sigma,\mu,\mu+0.5\sigma,\mu+\sigma\}=\{0.2,0.35,0.5,0.65,0.8\}$.

Note that in the look-ahead model both the expansion and production decisions must be made \textit{before} demand is realized in each stage, while distribution decisions are made \textit{after} demand is realized in each stage.  In practice, this means that the distribution decisions can always be optimally adjusted based on the availability of true information of the uncertainty. Further, a new set of first-stage (or ``here and now'') decisions can be obtained. This motivates us to primarily assess the quality of an optimal production policy using what we will refer to as ``Pseudo Simulator". In a pseudo simulator, the expansion decisions ($y^{\textrm{bin}}$) and the production decisions (i.e., the decision variables associated with the policies $x_{it}(\vxi_{[t]})$ and $x^{\rm{exp}}_{it}(\vxi_{[t]})$) are restricted to the optimal DR computed in the model, while allowing the distribution decisions (i.e., the decision variables associated with the policy $y_{ijt}(\vxi_{[t]})$) to be re-optimized within the PLDR-5[$\mathcal{Z}_{\textrm{base}}$] model. The pseudo simulator-based profit will be used to assess the flexibility of DRs. In our study, we assume that PLDR-5[$\mathcal{Z}_{\textrm{base}}$] model is equivalent to taking decisions given perfect information. 

\begin{figure}[H]
\begin{center}
\includegraphics[scale=0.5]{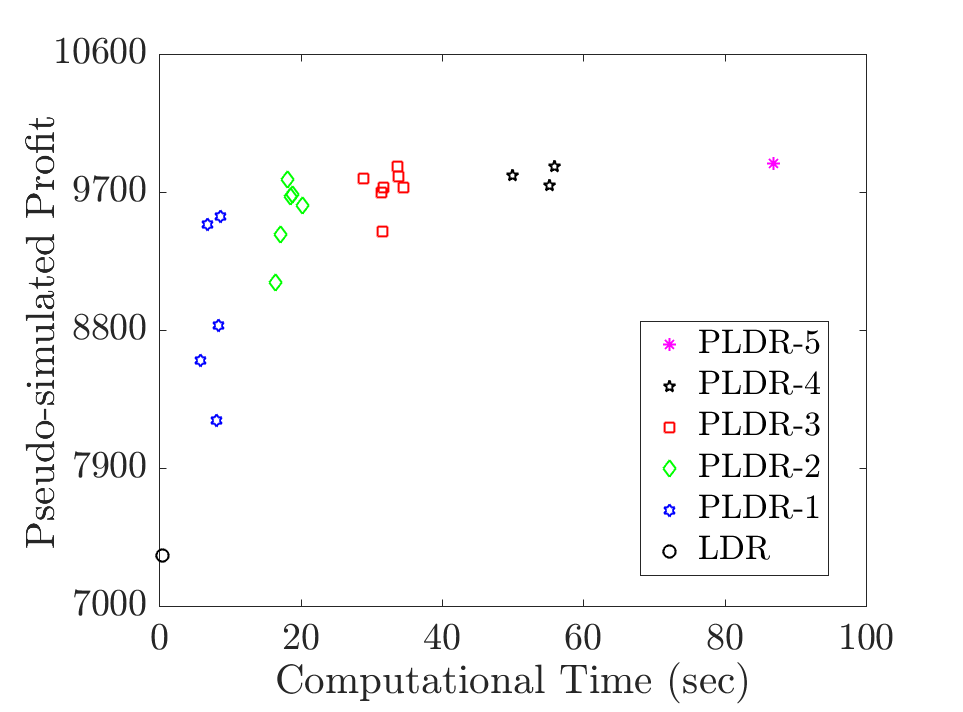}
\caption{Clusters of all possible PLDRs using 1, 2, 3, 4 and 5 breakpoints generated by solving the multistage stochastic transportation problem. The increase in the pseudo-simulated profit at the expense of the increased computational time is not attractive beyond PLDR-2, hence we will focus on designing competitive HDRs to LDR, PLDR-1 and PLDR-2. The solution quality of PLDR clusters with higher uncertainty resolution is less sensitive to the choice of the breakpoints. Parameters used: $T=10,\ \mathcal{Z}_{\textrm{base}}=\{0.2,0.35,0.5,0.65,0.8\},\ S_i= 0.15(C_i +H_i)\ \forall i \in \I,\ \forall \xi_t \sim \mathcal{U}(0,1)\ t \in \T_{-1}$.   }\label{fig:HDR_T10_PLDRs_clusters}
\end{center}
\end{figure}

\begin{figure}[H]
\begin{subfigure}{0.48\textwidth}
\centering
\includegraphics[width=7.5cm,height=6cm]{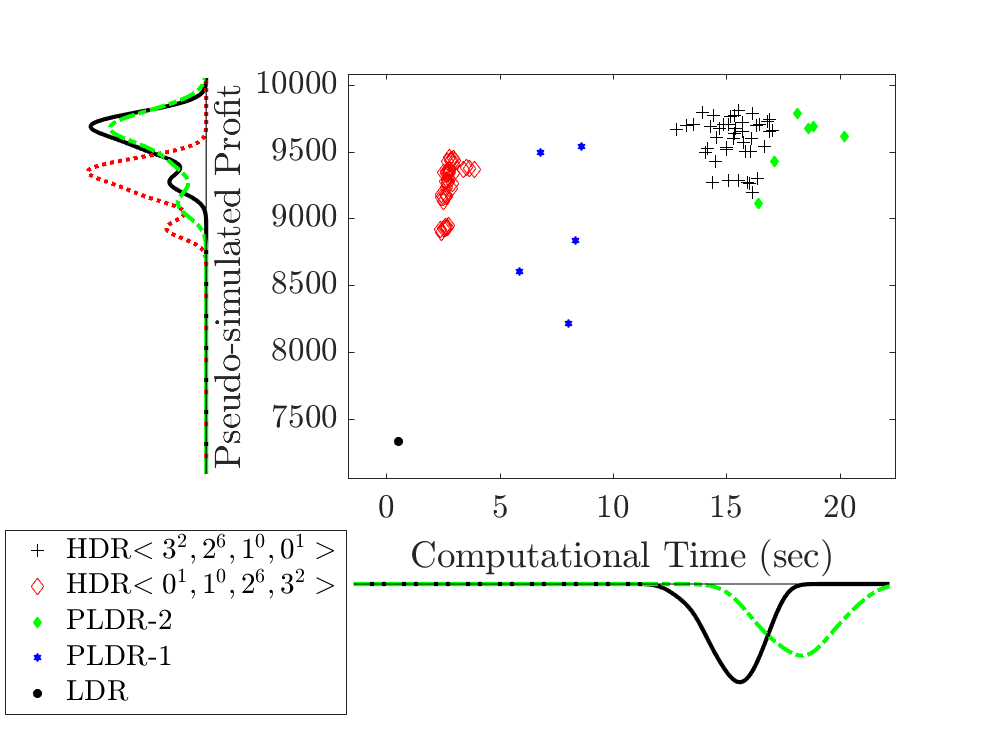}
\caption{}\label{fig:HDR_T10_a}
\end{subfigure}
\begin{subfigure}{0.48\textwidth}
\centering
\includegraphics[width=7.5cm,height=6cm]{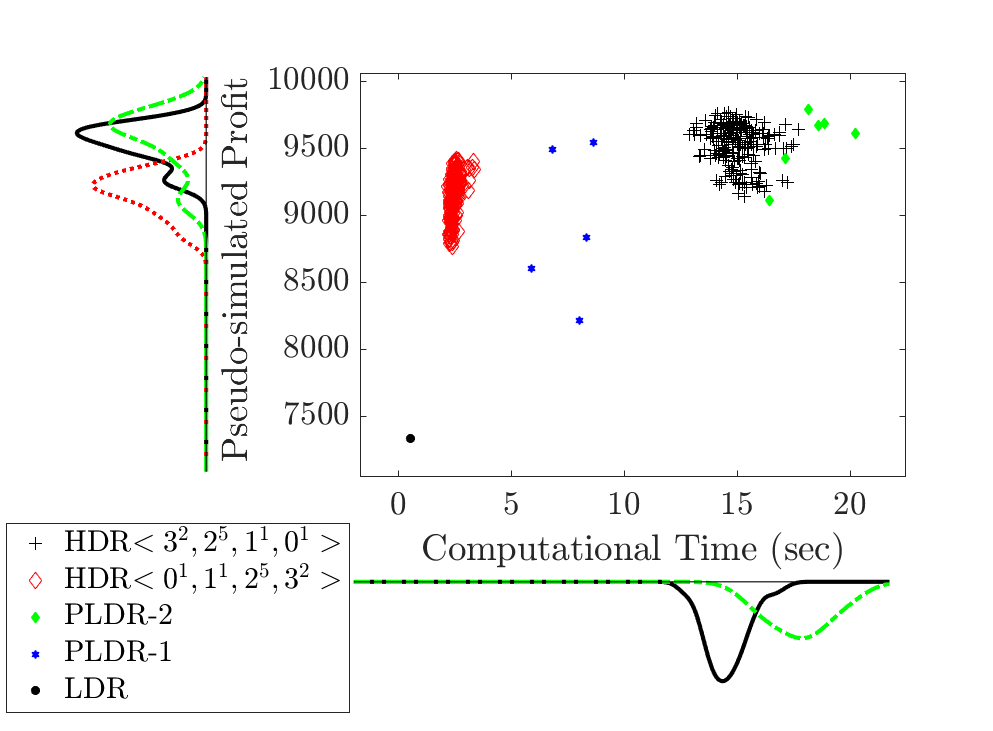}
\caption{}\label{fig:HDR_T10_b}
\end{subfigure}
\\
\begin{subfigure}{0.48\textwidth}
\centering
\includegraphics[width=7.5cm,height=6cm]{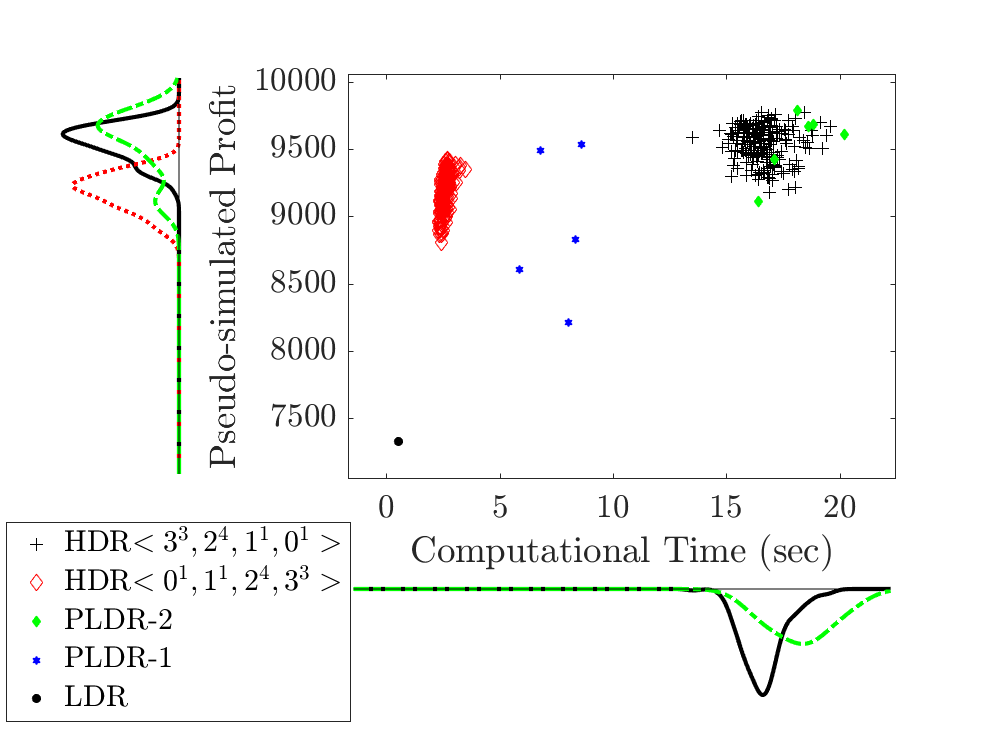}
\caption{}\label{fig:HDR_T10_c}
\end{subfigure}
\begin{subfigure}{0.48\textwidth}
\centering
\includegraphics[width=7.5cm,height=6cm]{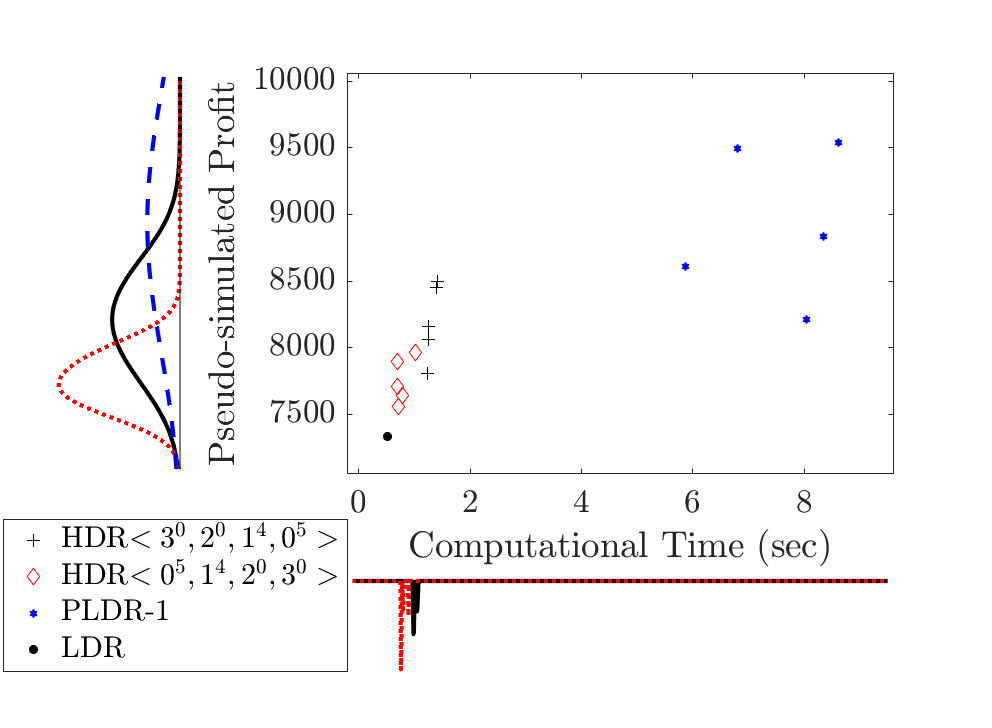}
\caption{}\label{fig:HDR_[111100000]}
\end{subfigure}
\caption{Pseudo-simulated profit vs computational time clusters are generated for the multistage stochastic transportation problem using non-increasing and non-decreasing HDRs. Non-increasing HDRs in Figs.~(a), (b) and (c) exhibit a competitive solution quality with respect to PLDR-2, whereas the deterioration of solution quality exhibited by non-decreasing HDRs is prominent. In Fig. (d), the non-increasing HDR is designed to improve the inferior LDR solution quality at the expense of a small increase in computational time. It is also more attractive than the non-decreasing HDR. Parameters used: $T=10,\ S_i= 0.15(C_i +H_i)\ \forall i \in \I,\ \xi_t \sim \mathcal{U}(0,1)\ \forall t \in \T_{-1}$. }
\label{fig:HDR_1}
\end{figure}

Given $\mathcal{Z}_{\textrm{base}}$, the total possible definitions of PLDR-1, PLDR-2, PLDR-3  and PLDR-4 are equal to 5, 10, 10 and 5, respectively. In order to design a comprehensive computational study, we impose some heuristics on the set of breakpoints used in a PLDR. In particular, the steps (i.e., $0.5\sigma$) between the minimum and maximum breakpoint values for PLDR-2, PLDR-3 and PLDR-4 are set to be at least 2, 3 and 4,  respectively. The computational time and pseudo-simulated profit for the LDR and PLDRs are shown in Figure~\ref{fig:HDR_T10_PLDRs_clusters}. Most of the improvement in the solution quality is achieved when the flexibility of the DR increases from an LDR to a PLDR-1. This comes at a relatively small increase in computational time. The increase in the solution quality with respect to the computational time post PLDR-2 is not very attractive. On a minor note, the clusters of PLDRs with higher flexibility (i.e., higher uncertainty resolution) are less dispersed, which indicates that the solution becomes less sensitive to the choice of the set of breakpoints. Stated differently, the deviation in solution quality from that of the true solution generated by a PLDR with higher number of breakpoints is less, given that (i) we do not have a method to select the optimal set of breakpoints and (ii) we sample them from $\mathcal{Z}_{\textrm{base}}$ instead.  
Defining a competitive non-increasing lifting strategy for an HDR is a design problem. One intuitive way to do so, is to use a PLDR as a template and simply increment and reduce the uncertainty resolution in early and late stages, respectively. Figure~\ref{fig:HDR_1} depicts the quality of four HDRs in a pseudo-simulated profit vs computational time plot. The non-increasing HDRs are compared with the relative non-decreasing HDRs (i.e., the inverse). For notation purposes, HDR${<3^2,2^6,1^0,0^1>}$ indicates that 3 breakpoints are used to lift $\xi_t$ in the first two stages, then 2 breakpoints are used to lift $\xi_t$ in the next six stages, none of $\xi_t$ is lifted with 1 breakpoint and $\xi_{T}$ in the last stage is not lifted. In the other hand HDR${<0^1,1^0,2^6,3^2>}$ corresponds to the inverse lifting strategy.

  The computational experiments are comprehensive and include all the possible sets of breakpoints for each lifting resolution in an HDR (i.e., 42 combinations for HDR$<3^2,2^6,1^0,0^1>$). The solution quality exhibited by  HDR $<3^2,2^6,1^0,0^1>$, HDR $<3^2,2^5,1^1,0^1>$ and HDR $<3^3,2^4$ $,1^1,0^1>$ are competitive with PLDR-2 at a reduced computational time. Despite the higher reduction in computational time generated by the non-decreasing HDRs, the deterioration in solution quality is significant which makes them not attractive with respect to PLDR-2. In Figure~\ref{fig:HDR_1}(d), HDR$<3^0,2^0,1^4,0^5>$ elevates the inferior solution quality of LDR at a minimal additional computational cost. Such HDRs represent a competitive alternative to PLDR-1 for problems where the latter DR is computationally prohibitive. 

\begin{figure}[H]
\begin{subfigure}{0.48\textwidth}
\centering
\includegraphics[width=7.5cm,height=6cm]{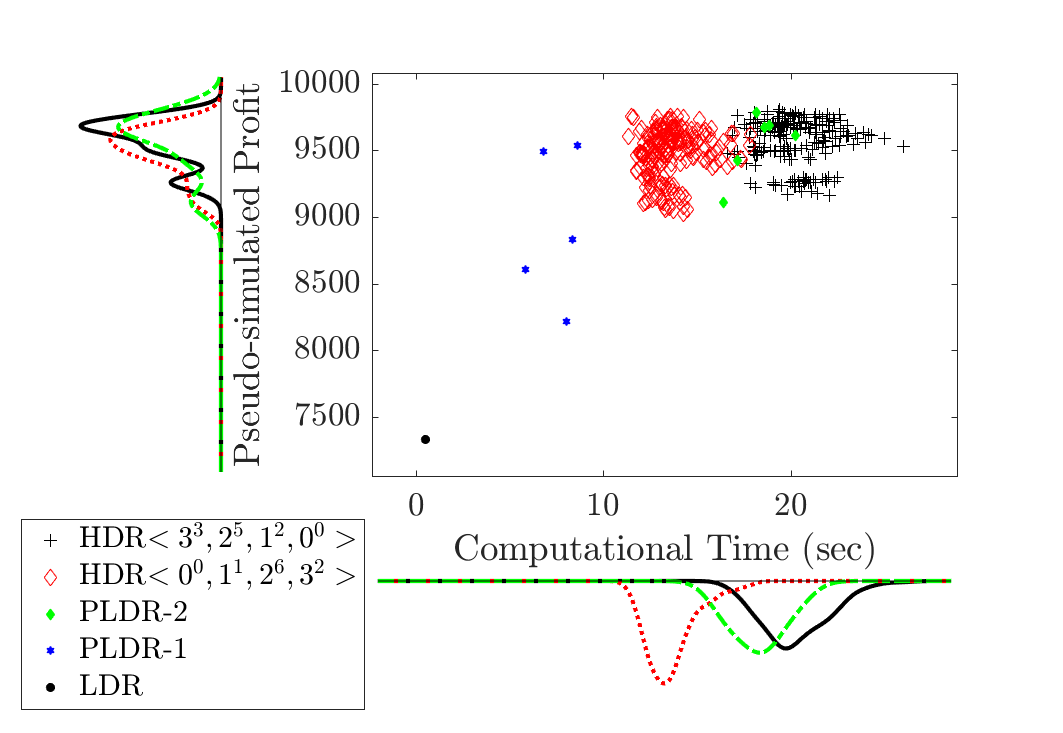}
\end{subfigure}
\begin{subfigure}{0.48\textwidth}
\centering
\includegraphics[width=7.5cm,height=6cm]{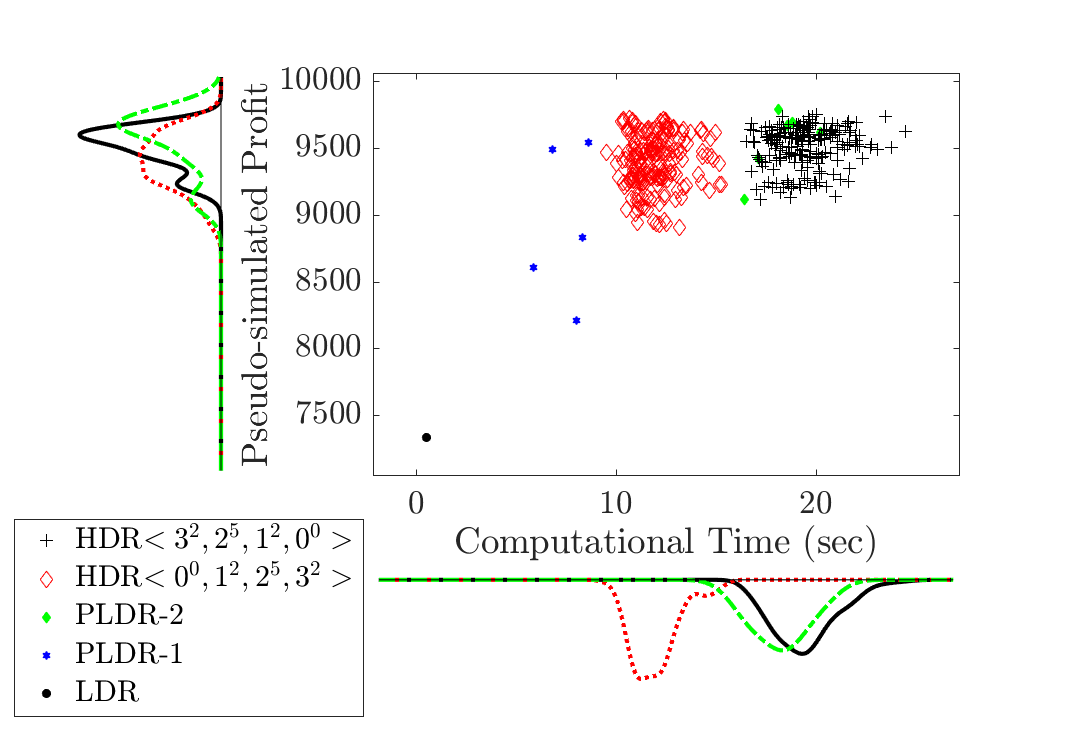}
\end{subfigure}
\\
\begin{subfigure}{0.48\textwidth}
\centering
\includegraphics[width=7.5cm,height=6cm]{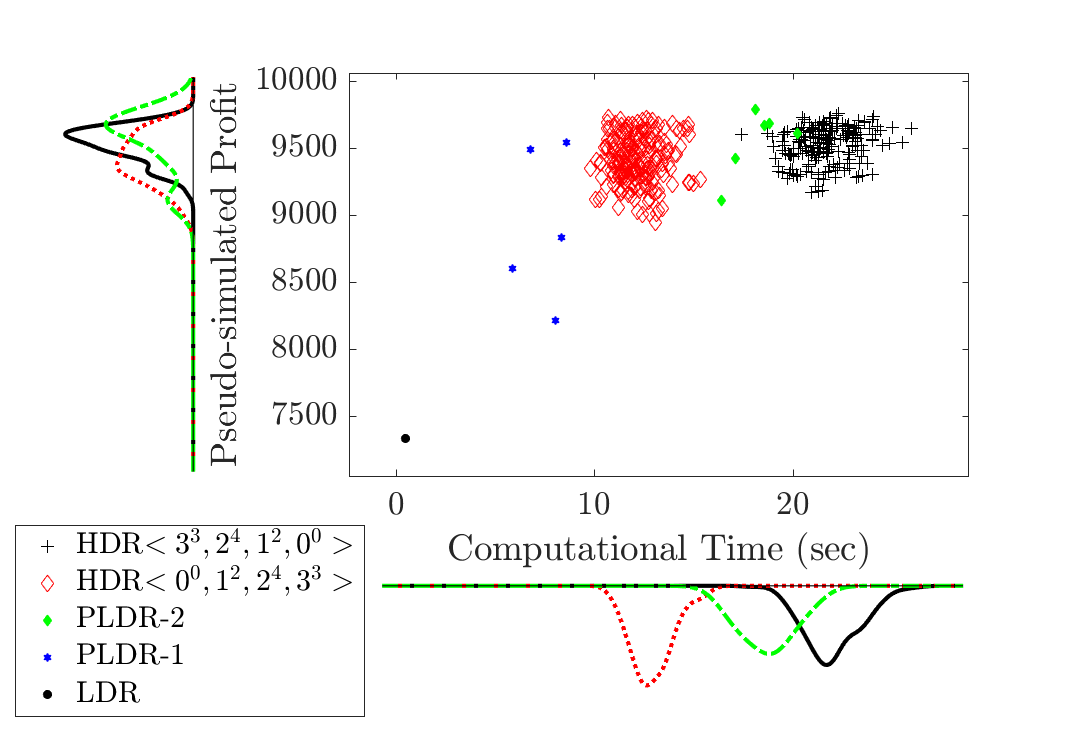}
\end{subfigure}
\caption{Pseudo-simulated profit vs computational time clusters are generated for the multistage stochastic transportation problem using poorly designed non-increasing HDRs. The HDRs are similar to those in Figs.~\ref{fig:HDR_1}(a)-\ref{fig:HDR_1}(c), but with a resolution of 1 instead of 0 in the last stage. It is clear that the additional computational burden is not met with a corresponding increase in solution quality. In fact, the non-increasing HDRs lose their competitive advantage with respect to PLDR-2. Further, it is shown that the non-decreasing HDRs offer a better trade-off between between solution quality and computational time, nevertheless this is when compared to a poorly designed non-increasing HDR. Parameters used: $T=10,\ S_i= 0.15(C_i +H_i)\ \forall i \in \I,\ \xi_t \sim \mathcal{U}(0,1)\ \forall t \in \T_{-1}$.}\label{fig:HDR2_T10}

\end{figure}

The quality of an HDR is dictated by the design of the lifting strategy. A poorly designed non-increasing HDR will exhibit unsatisfactory performance. This is observed in Figure~\ref{fig:HDR2_T10} where we implement the same HDRs in Figs.~\ref{fig:HDR_1}(a)-\ref{fig:HDR_1}(c), but with a trailing (i.e., in the last stage) resolution of 1 instead of 0. A wasted increase in the computational time that is not accompanied with an increase in solution quality is best seen for  HDR${<3^3,2^5,1^2,0^0>}$. This supports the finding that lifting in late stages is not recommended as it adds to the computational burden without any significant quality return. On the other hand, the non-decreasing HDRs now have a leading resolution of 1, which is enough to describe the uncertainty in early stages at a sufficient level of detail (i.e., see Figure \eqref{fig:HDR_T10_PLDRs_clusters}). Consequently, the solution quality gap with non-increasing HDRs is significantly reduced. 

\subsubsection*{Implementing HDRs for a 20-stage planning horizon}
\begin{figure}[H]
\begin{center}
\includegraphics[scale=0.5]{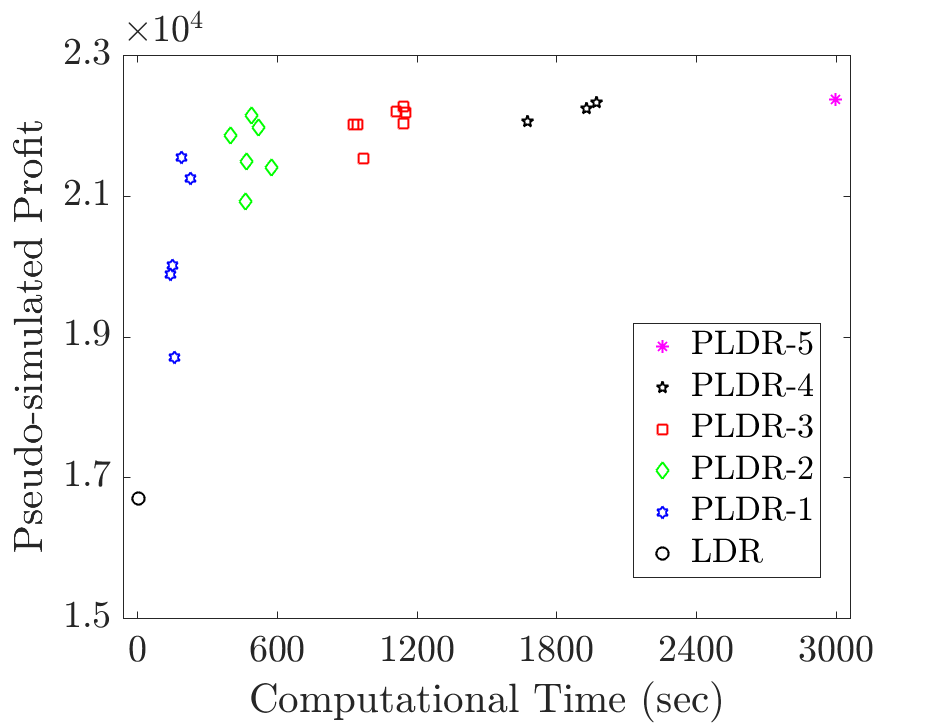}
\caption{Pseudo-simulated profit vs computational time clusters of PLDRs-$i$, $i \in \{1,2,3,4\}$ generated by solving a multistage stochastic transportation problem for $T=20$. The increase in the profit with respect to computational time is not attractive beyond PLDR-2, hence we will focus on designing competitive non-increasing HDRs with respect to LDR, PLDR-1 and PLDR-2. Parameters used: $ S_i= 0.7(C_i +H_i)\ i \in \I,\ \xi_t \sim \mathcal{U}(0,1)\ t \in \T_{-1}$.   }\label{fig:HDR_T20_PLDRs_clusters}
\end{center}
\end{figure}

In this section, we show the applicability of non-increasing HDRs in settings with a more intensive computational burden. To do so, we extend the planning horizon $T$ to 20, and we keep the rest of the computational setting unchanged. Figure~\ref{fig:HDR_T20_PLDRs_clusters} shows the LDR and PLDR solution clusters  in terms of pseudo-simulated profit and computational time. There is no change to the heuristics used to select the breakpoints for PLDRs. The computational time has increased by at least one order of magnitude, compared to the 10-stage problem. The improvement in pseudo-simulated profit beyond PLDR-2 is minimal, thus we instead explore competitive HDRs to LDRs, PLDRs-1 and PLDRs-2.
With a larger planning horizon, the design of a competitive non-increasing HDR becomes more intricate. Over-lifting the uncertainty in early stages will accumulate unnecessary computational overhead (i.e., with no solution quality return). To design an HDR, we first identified the best set of breakpoints corresponding to each PLDR-$i$, $i=\{1,2,3,4\}$ (see Table \ref{table:cost_comparison_hybrid_PLDR_T10} in appendix \ref{appendix_B}). Then, we investigate the impact of the lifting strategy on the trade-off between solution quality and computational time through four computational experiments.

\begin{figure}[H]
\begin{subfigure}{0.48\textwidth}
\centering
\includegraphics[width=7cm,height=5cm]{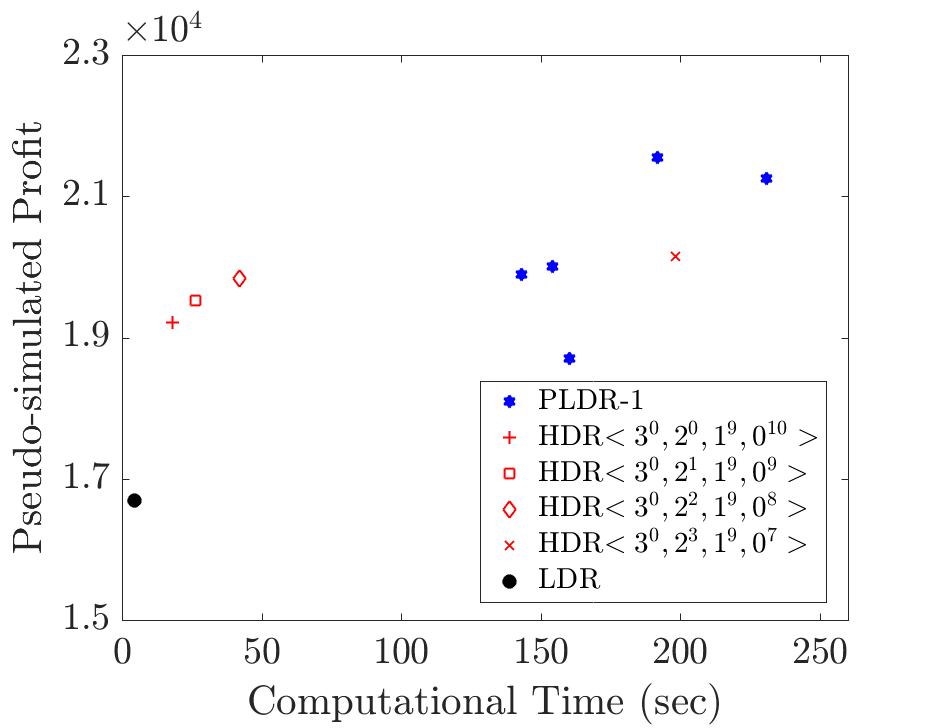}
\caption{}\label{fig:HDR_T20_a}
\end{subfigure}
\begin{subfigure}{0.48\textwidth}
\centering
\includegraphics[width=7cm,height=5cm]{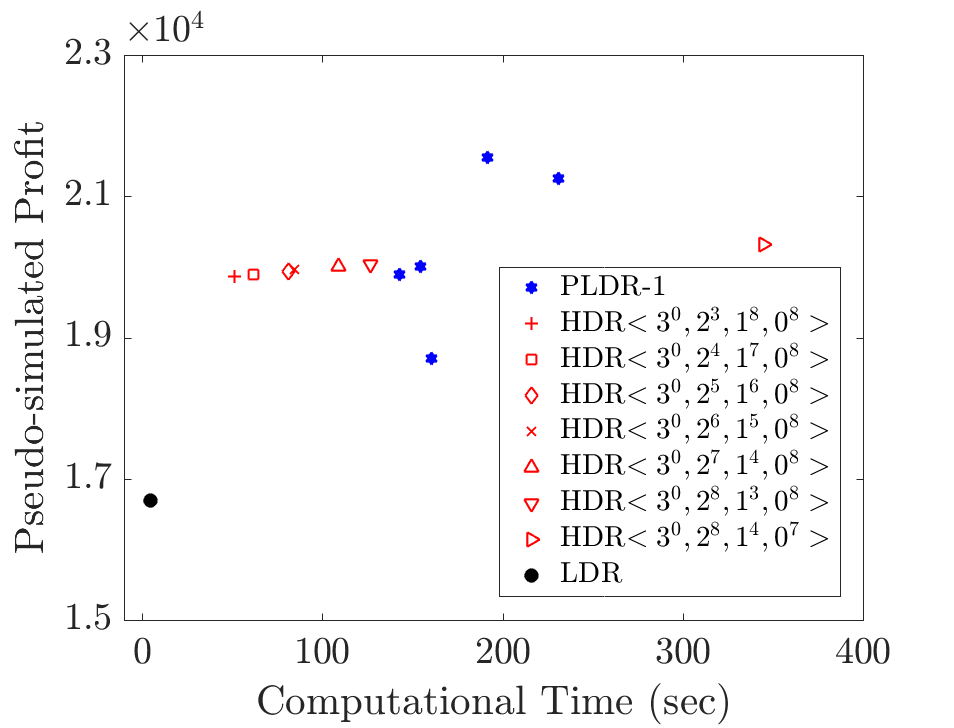}
\caption{}\label{fig:HDR_T20_b}
\end{subfigure}\\
\begin{subfigure}{0.48\textwidth}
\centering
\includegraphics[width=7cm,height=5cm]{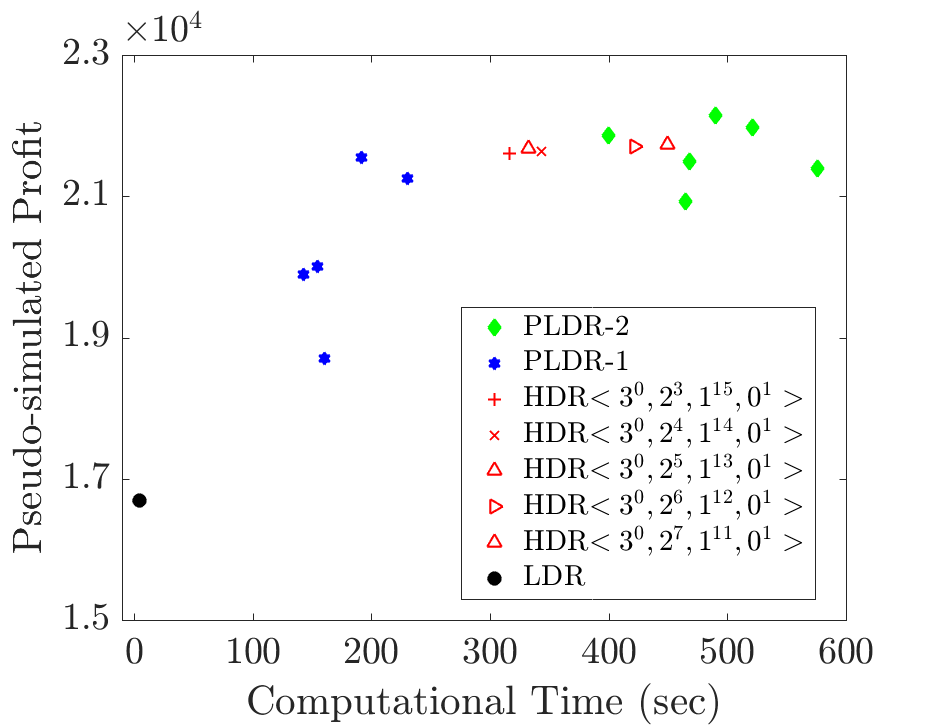}
\caption{}\label{fig:HDR_T20_c}
\end{subfigure}
\begin{subfigure}{0.48\textwidth}
\centering
\includegraphics[width=7cm,height=5cm]{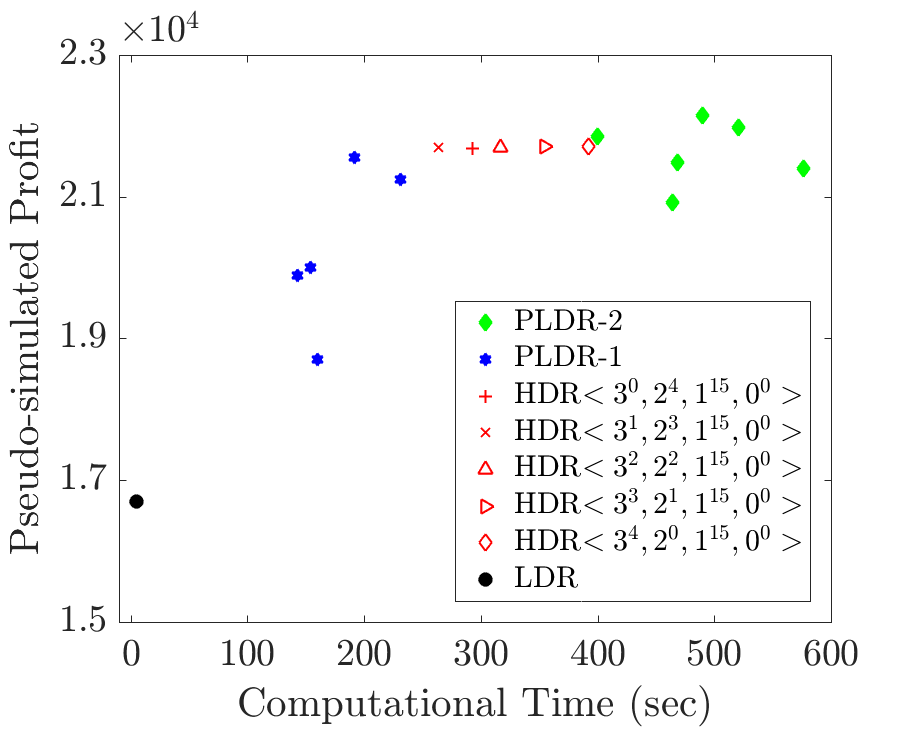}
\caption{}\label{fig:HDR_T20_d}
\end{subfigure}

\caption{Pseudo-simulated profits vs computational time using four sets of non-increasing HDRs in a multistage stochastic transportation problem. For all sets, an HDR with a good lifting strategy design offers computational benefits in terms of the trade-off between solution quality and computational time. However, as in any design problem, a poorly designed non-increasing HDR loses the aforementioned competitive advantage. Parameters used: $T=20,\ S_i= 0.7(C_i +H_i)\ \forall i \in \I,\ \xi_t \sim \mathcal{U}(0,1)\ \forall t \in \T_{-1}$.}
\label{fig:HDR_T20}
\end{figure}

Figure \ref{fig:HDR_T20} depicts the results of four sets of non-increasing HDRs with different lifting strategies. Figures \ref{fig:HDR_T20}(a) and \ref{fig:HDR_T20}(b) showcase systematic search for competitive non-increasing HDRs that improve the inferior solution quality of an LDR at the expense of a small increase in the computational time. Still, a poorly designed HDR does lose the competitive advantage it may exhibit as in the case of HDR${<3^0,2^3,1^9,0^7>}$ and HDR${<3^0,2^8,1^4,0^7>}$. Similarly, in Figures \ref{fig:HDR_T20}(c) and \ref{fig:HDR_T20}(d), two sets of non-increasing HDRs are defined to capture a portion of the increase in the solution quality between PLDR-1 and PLDR-2 clusters, at the expense of a partial increase in the computational time. Based on the design of the lifting strategy, most of the HDRs are deemed attractive except for a couple HDRs in Fig.~\ref{fig:HDR_T20}(c) and an HDR in Fig.~\ref{fig:HDR_T20}(d) which are poorly designed.

\section{Conclusion}
Since their recent inception, piecewise linear decision rules have received increasing attention in the stochastic and robust optimization communities.  This is mainly due to the flexibility and solution quality improvements that they provide, while maintaining a tractable linear structure closely resembling that used for LDRs. However, the increase in solution quality comes at the expense of a significant computational burden for large-scale practical problems. This study first provides an unexpected result in which an LDR is \textit{superior} to a PLDR with a single breakpoint when assessed within a simulation environment. The comparison is done using a hyper-rectangle uncertainty set where there is an exact tractable representation of the convex hull of its lifted uncertainty set (i.e., no overestimation). This finding highlights the need for assessing the quality of the policies computed by the look-ahead model within a simulator environment, instead of just relying on the objective function value. Then, the study emphasizes on the concept of implementing hybrid decision rules as a promising direction to mitigate the increased computational burden in practical multistage adaptive optimization problems. The main aspect of HDRs explored is the lifting strategies or the axial combinations of the LDR and PLDRs where it is empirically illustrated that having higher uncertainty resolution (i.e., more linear pieces) in early stages is more important than having it in later stages.
 
There are several open questions that are ripe for additional research. First, the flexibility of a non-increasing HDR highly depends on the design of the lifting strategy. Currently, we do not provide a systematic way to design an optimal lifting strategy. Second, this computational study addresses only one-dimensional parameter uncertainty in each stage. The computational benefits manifested by HDRs can be obtained at a higher dimension of $\vxi_t$, but the design of a good lifting strategy becomes more challenging. Third, removing the assumptions on the general model will increase the complexity of the problem. Cost coefficients uncertainty lead to quadratic terms in the objective function. An uncertain recourse matrix generates nonlinear uncertain constraints whose tractable stochastic/robust counterparts are more challenging to derive. As for defining binary variables as adaptive decisions, it has seen limited applications (see, e.g., \cite{bertsimas2015design}) in the context of decision rule-based methods. Lastly, the hybrid decision rule approach can be extended to conditional value at risk and robust multistage adaptive optimization problems.
 
\subsection*{Acknowledgement}
The authors would like to thank Nikolaos Lappas for his comments and feedback. 

\section*{Nomenclature}
\subsection*{General}
\begin{tabular}{ll}
\toprule
& Definitions \\
\hline
\textbf{Abbreviations}& \\
DR   & Decision rule\\
LDR  & Linear DR\\
PLDR-$i$ ($\vz$) & Piecewise LDR with $i$ breakpoints at $\vz$ in each stage\\
HDR${<3^l,2^k,1^p,0^s>}$   & Hybrid DR where the uncertainty in the first $l$ stages is lifted by $3$ \\ 
  & breakpoints; the next $k$ stages   is lifted by $2$ breakpoints; the next $p$   \\
  & stages is lifted by $1$ breakpoint, and the last $s$ stages is not lifted \\
LASC  & Linear adaptive stochastic counterpart \\
PWLASC  & Piecewise Linear adaptive stochastic counterpart \\
conv  & convex hull operator\\
\\
\textbf{Sets} &\\
$\T$   & Time stages \\
$\T_{-t}$ & Time stages excluding stage $t$ \\
$\Xi$ & Uncertainty set\\
$\Xi_t^\prime$ & Lifted uncertainty in stage $t$  \\
$\Xi^\prime$ & Overall Lifted uncertainty set   \\
$\bar{\Xi}^\prime$ & Outer approximation of conv $\Xi^\prime$ \\
$\hat{\Xi}_t$ & Scenario-based uncertainty set in stage $t$   \\
\\
\textbf{Parameters} &\\
$\xi_t$ & Primitive uncertainty in stage $t$ \\
$\hat{\xi}_t$ & Scenarios approximating $\xi_t$   \\
$\vW$ & Matrix of $\Xi$ \\
$\vh$ & Right hand side vector of $\Xi^\prime$ \\
$\vA_t$ & Matrix of conv  $\Xi^\prime_t$ \\
$\vb_t$ & Right-hand side vector of conv $\Xi^\prime_t$ \\
$\vA^{\rm{l}}$ & Matrix of $\bar{\Xi}^\prime$ \\
$\vb^{\rm{l}}$ & Right-hand side vector of $\bar{\Xi}^\prime$ \\
$\vw_i$ & The $i^{th}$ column vector of $\vW$ \\
$\vV_t$ & Observation matrix in stage $t$ \\
$\ve_{t}$  &  Vector with a value $1$ at the $t$ index, and 0 otherwise\\
$\ve^\prime_{t}$ & Vector with a value of 1 from the index  $\sum_1^{t-1}r_t +1$ to $\sum_1^{t}r_t$,\\
 & and 0 otherwise  \\    
$\sigma$ & Standard deviation  \\   
\bottomrule
\end{tabular}

Unless otherwise stated, the superscript `` $^\prime$ "  refers to the same variable definition, but in the lifted uncertainty space.
\subsection*{Newsvendor Model}
\begin{tabular}{ll}
\toprule
 & Definitions \\
 \hline
\textbf{Variables} \\
$x_t$ & Ordered units in stage $t$ \\
$I_t$      & Balance of units in stage $t$\\
$s_t^+$    & Inventory units in stage $t$ \\
$s_t^-$      &  Backlog units  in stage $t$ \\
$p_t^0$ & Intercept of LDR in stage $t$ where $p=\{x,I,s^+,s^-\}$\\
$\v{P}_t^1$ & Slope of  LDR at stage $t$ where $\v{P}=\{\vX,\vI,\vS^+,\vS^-\}$ \\
\\
\textbf{Parameters} & \\
$d_t$   & Demand in stage $t$ \\
$r_t$   &  Number of breakpoints for $d_t$ \\
$z_j^i$ & The $j^{th}$ breakpoint  in $d_i$ \\
$d_{tj}^\prime$ &  The $j^{th}$ lifted  element of $d_t$  \\
$U^x$ & Ordering amount limit \\
$C_t$ & Purchasing cost in stage $t$ \\
$H_t$ & Inventory cost in stage $t$ \\
$B_t$ & Backlogging cost in stage $t$ \\
$I_1$ & Initial inventory \\
$l_t,\ u_t$ & Lower and upper bound of $d_t$ \\
\bottomrule
\end{tabular}

\subsection*{Transportation Model}
\begin{tabular}{ll}
\toprule
 & Definitions \\
\hline
\textbf{Sets} & \\
$\I$   & Set of suppliers   \\
$\J$   & Set of customers \\
$\mathcal{Z}_{\textrm{base}}$ & Base set of potential breakpoints \\ 
\\
\textbf{Variables} & \\
$x_{it}$   &  Produced units by supplier $i$ in stage $t$ \\
$x^{\textrm{exp}}_{it}$  &  Produced units within the expanded capacity by supplier $i$ in stage $t$ \\
$I_{it}$   & Inventory stored by supplier $i$ in stage $t$\\
$y_{ijt}$  & Transported items from supplier $i$ to customer $j$ in stage $t$ \\
$p_t^0$    & Intercept of LDR at stage $t$ where $p:=\{x_i,I_i,y_{ij}\}$ \\
$\v{P}_t^1$ & Slope of  LDR at stage $t$ where $\v{P}:=\{\vX_i,\vI_i,\vy_{ij}\}$ \\
$y^{\textrm{bin}}_{i}$ & Equals 1 when capacity is expanded in supplier $i$; otherwise, 0\\
          \\
\textbf{Parameters} & \\

$R_{jt}$   & Unit revenue of customer $j$ in stage $t$ \\
$T_{ijt}$  & Unit transportation cost along $(i,j)$ arc in stage $t$\\
$C_{it}$   &  Unit production cost of supplier $i$ in stage $t$ \\
$M_{i}$    & Capital expansion cost of supplier $i$  \\
$Q_{i}$    & Additional capacity for supplier $i$  \\
$H_{it}$   & Unit holding cost of supplier $i$ in stage $t$ \\
$S_i$      & Unit salvage value for supplier $i$ \\
$U_i^{\max}$ & Production limit of supplier $i$\\
$D_{jt}$   & Customer demand $j$ in stage $t$\\
$D^0_{jt}$ & Intercept of customers $j$'s linear demand function in stage $t$\\
$D^1_{jt}$ & Slope of customers $j$'s linear demand function in stage $t$\\
\bottomrule
\end{tabular}

{\fontsize{10}{12}{\selectfont
\bibliographystyle{abbrvnat}

}}

\newpage
\appendix
\section{Policies of multistage stochastic newsvendor problem}
The ordering, inventory and backlog policies for the multistage stochastic newsvendor problem in section~\ref{sec:newsvendor_numerical_results} are depicted in the following table.
{\fontsize{9}{11}\selectfont { 
  \begin{table}[H]
  \centering
  \caption{Optimal adpative policies for the multistage stochastic newsvendor problem using an LDR, a PLDR-1 ($\E[d_t]$) and a PLDR-1 ($U^x$). Parameters used: $T=4,\ U^x=8,\ I_1=4, d_t\sim \mathcal{U}(0,1)\ t \in \T_{-1}$. } \label{tab:optimal_policies_multistage_inventory_solution_T_4}
    \begin{tabular}{c}
    \toprule
    $\begin{bmatrix} x_1 \\  x_2 \\ x_3 \\ s^+_2 \\  s^+_3 \\ s^+_4 \\  s^-_2 \\ s^-_3 \\ s^-_4 \end{bmatrix} 
             = \begin{bmatrix}0 & 0 & 0 \\ 0.8& 0& 0 \\ 0& 0.8& 0 \\ -1 & 0 & 0 \\ -0.2& -1& 0 \\ 0& -0.2& -1\\  0 & 0 & 0 \\ 0 & 0 & 0 \\ 0.2& 0& 0\end{bmatrix} \begin{bmatrix}
             d_2\\ d_3\\ d_4 \end{bmatrix} + \begin{bmatrix} 8\\ 0\\ 0 \\ 12 \\ 12 \\ 12 \\ 0\\ 0\\ 0 \end{bmatrix}$   
         \\
            (a) LDR \\
          $\begin{bmatrix} x_1 \\  x_2 \\ x_3 \\ s^+_2 \\  s^+_3 \\ s^+_4 \\  s^-_2 \\ s^-_3 \\ s^-_4 \end{bmatrix} 
             = \begin{bmatrix}0 & 0 & 0&0 & 0 & 0 \\ 0.6& 1& 0& 0& 0& 0 \\ 0& 0& 0.4& 1& 0& 0 \\-1 & -1 & 0 & 0& 0 & 0 \\ -0.4& 0& -1& -0.6 &0 &0 \\ -0.4 &0 &-0.6 &0 &-1 &0 \\ 0 & 0 & 0& 0 & 0 & 0 \\ 0 & 0 & 0& 0.4 & 0 & 0 \\ 0 & 0 & 0& 0 & 0 & 1 \end{bmatrix} \begin{bmatrix}
             d^\prime_{21}\\ d^\prime_{22}\\ d^\prime_{31}\\d^\prime_{32}\\d^\prime_{41}\\d^\prime_{42} \end{bmatrix} + \begin{bmatrix} 6\\ 0\\ 0 \\ 10 \\ 10 \\ 10 \\ 0\\ 0\\ 0 \end{bmatrix}$ 
           \\
    (b) PLDR-1 ($\E[d_t]$) \\
              $\begin{bmatrix} x_1 \\  x_2 \\ x_3 \\ s^+_2 \\  s^+_3 \\ s^+_4 \\  s^-_2 \\ s^-_3 \\ s^-_4 \end{bmatrix} 
             = \begin{bmatrix}0 & 0 & 0&0 & 0 & 0 \\ 1& 0& 0& 0& 0& 0 \\ 0& 0& 1& 0& 0& 0 \\-1 & 0 & 0 & 0& 0 & 0 \\ 0& 0& -1& 0 &0 &0 \\ 0 &0 &0 &0 &-1 &0 \\ 0 & 1 & 0& 0 & 0 & 0 \\ 0 & 1 & 0& 1 & 0 & 0 \\ 0 & 1 & 0& 1 & 0 & 1 \end{bmatrix} \begin{bmatrix}
             d^\prime_{21}\\ d^\prime_{22}\\ d^\prime_{31}\\d^\prime_{32}\\d^\prime_{41}\\d^\prime_{42} \end{bmatrix} + \begin{bmatrix} 4\\ 0\\ 0 \\ 8 \\ 8 \\ 8 \\ 0\\ 0\\ 0 \end{bmatrix}$ 
           \\
    (c) PLDR-1 ($U^x$) \\
    \bottomrule
    \end{tabular}
  \end{table}}
}
\section{Computational setting for multistage transportation problem with expansion decisions}\label{appendix_B}
The computational parameters used for the multistage stochastic transportation problem with expansion decisions in section~\ref{sec:transporttion_numerical_results} are shown in the following two tables.
\begin{table}[H]
\caption{Computational parameters for the multistage stochastic transportation problem with 10 suppliers and 10 customers. They are assumed to be constant for all stages (e.g. $C_{it} \equiv C_i\ \forall t$).}
\begin{center}
\begin{tabular}{lccccccccccccccc}
\toprule
  & $C_{i}$ & $H_{i}$ & $M_i$ & $U_i^{\max}$ & $Q_i$ & \multicolumn{10}{c}{$T_{ij}$}  \\
\cmidrule(lr){7-16} 
$i \downarrow j \rightarrow$ & & & & & & 1 & 2 & 3 & 4 & 5 & 6 & 7 & 8 & 9 & 10  \\
\hline
 1   & 4    & 1   & 6    & 10& 3 & 2 & 9 & 11 & 5& 8& 13&6&9&12&7   \\
 2   & 7    & 2   & 12   & 8 & 7 & 5&1&6&8&12&10&9&11&5&8  \\
 3   & 6.50 & 2   & 10   & 5 & 5 & 9&13&2&11&7&6&13&7&10&13   \\
 4   & 3    & 0.75& 8    & 6 & 4& 13&5&9&3&9&12&7&10&13&6 \\
 5   & 1    & 0.50& 5    & 7 & 3& 7&10&13&6&1&9&10&12&6&10   \\
 6   & 8    & 3   & 14.50& 3 & 8& 10&11&5&9&13&3&11&5&8&12   \\
 7   & 5    & 1.50& 9    & 12& 4& 1&7&8&12&5&5&2&8&11&5   \\
 8   & 4.25 & 1   & 11   & 2 & 5& 6&8&12&7&10&11&5&1&7&9   \\
 9   & 6    & 2   & 8    & 4 & 4& 8&12&7&10&11&8&8&13&3&11 \\
 10  & 2    & 0.50& 5    & 4 & 9& 12&6&10&12&6&7&12&6&10&2\\
\bottomrule
\end{tabular}
\end{center}
\end{table}

\begin{table}[H]
\caption{Demand and revenue parameters for the multistage stochastic transportation problem with 10 suppliers and 10 customers. They are assumed to be constant for all stages (e.g. $D^0_{jt} \equiv D^0_j\ \forall t$).}
\begin{center}
\begin{tabular}{lccccccccccccc}
\toprule
 $j$  & $D^0_{j}$ & $D^1_{j}$  & $R_{j}$  \\
\hline
 1   & 6    & 3   & 17       \\
 2   & 5    & 2   & 23      \\
 3   & 2    & 2   & 15     \\
 4   & 7    & 1& 18      \\
 5   & 6.5  & 2.5& 19      \\
 6   & 4    & 2   & 16.5   \\
 7   & 8    & -1& 18       \\
 8   & 10   & -3   & 20     \\
 9   & 9    & 0.5   & 21      \\
 10  & 3    & 2& 12     \\
\bottomrule
\end{tabular}
\end{center}
\end{table}

\section{Best set of breakpoints for PLDRs in multistage transportation problem with 20 stages}
The following table illustrates the best set of breakpoints given $\mathcal{Z}_{\rm{base}}$ for the PLDRs investigated in the multistage stochastic transportation problem with expansion decisions. The values of the breakpoints were used to defined the uncertainty resolution in the HDRs explored for $T=20$ in section~\ref{sec:transporttion_numerical_results}.
\begin{table}[H]
\begin{center}
\caption{Best set of breakpoints obtained for the PLDRs investigated in a multistage stochastic transportation problem with expansion decisions. The main increase in solution quality occurs when the flexibility of the DR is increased from LDR to PLDR-1[0.65]. The increase in profit past PLDR-2[0.35,0.65] is less then $1\%$, hence it is not practical to design competitive HDRs beyond that. Parameters used: $T=20,\ \mathcal{Z}^{\rm{base}}=\{0.2,0.35,0.5,0.65,0.8\}, \ S_i= 0.7(C_i +H_i)\ \forall i \in \I,\ \xi_t \sim \mathcal{U}(0,1)\ \forall t \in \T_{-1}$.}
\label{table:cost_comparison_hybrid_PLDR_T10}
\begin{tabular}{lcccccccc}
\toprule
&     \multicolumn{2}{c}{First stage costs}    \\
\cline{2-3}
DRs  & Prod. & Exp.& Model &Pseudo Simulator  & Gap($\%)$ & Time (sec)\\\hline
LDR                          & 1449.50 & 675  & 15930  & 16692  & ---  & 4.42\\
PLDR-1[0.65]                 & 1298.40  & 518  & 21443    & 21550  & 29.10  & 191.86\\
PLDR-2[0.35,0.65]            & 1298.20  & 518  & 22054 & 22144 & 32.66 & 490.32\\
PLDR-3[0.2,0.5,0.65]         & 1298.50  & 518 & 22185 & 22271 & 33.42 & 1142.7\\                                                 PLDR-4[0.2,0.5,0.65,0.8]     & 1298.20 & 518 & 22323 & 22324 & 33.74 &1972.4 \\                                              
PLDR-5[$\mathcal{Z}^{\rm{base}}$]& 1297.5 & 518  & 22368 & --- & 34.00& 2997.6\\
                     
\bottomrule 
\end{tabular}
\label{tab:optm-dec-sol-T10}
\end{center}
\end{table}

\end{document}